\documentclass[11pt]{article}
\usepackage{amsmath,amsthm,amscd,amsfonts,amssymb}
\usepackage{graphicx}

\newtheorem{theorem}{Theorem}[section]
\newtheorem{lemma}[theorem]{Lemma}
\newtheorem{corollary}[theorem]{Corollary}
\newtheorem{proposition}[theorem]{Proposition}

\theoremstyle{definition}

\newtheorem{definition}[theorem]{Definition}
\newtheorem{example}[theorem]{Example}
\newtheorem{remark}[theorem]{Remark}

\newcommand{\const}{\mbox{const}}

\newcommand{\cl}{{\rm cl}}

\newcommand{\Log}{{\rm Log}}
\newcommand{\dist}{\mbox{\rm{dist}}}

\newcommand{\Imag}{{\rm Im}}
\newcommand{\Real}{{\rm Re}}
\newcommand{\speca}{{\rm spec}}
\newcommand{\Arg}{{\rm Arg}}
\newcommand{\Aut}{{\rm Aut}}
\newcommand{\rank}{{\rm rank}}

\begin{document}

\title{Holomorphic Semi-almost Periodic Functions}

\author{Alexander Brudnyi\thanks{Research supported in part by NSERC.
}\\
Department of Mathematics and Statistics\\
University of Calgary\\
Calgary, Alberta, T2N 1N4, CANADA\\
albru@math.ucalgary.ca\\
\and Damir Kinzebulatov\\
Department of Mathematics\\
University of Toronto\\
Toronto, Ontario, M5S 2E4, CANADA\\
dkinz@math.toronto.edu}

\date{}
\maketitle

%
%
%
%
%

\begin{abstract}
We study the Banach algebras of bounded holomorphic functions on the unit disk whose boundary values, having, in a sense, the weakest possible discontinuities, belong to the algebra of semi-almost periodic functions on the unit circle. The latter algebra contains as a special case an algebra introduced by Sarason 
in connection with some problems in the theory of Toeplitz operators. 

\end{abstract}

\section{Introduction}
\label{introsect}

We study the Banach algebras of \textit{holomorphic semi-almost periodic functions}, i.e., bounded holomorphic functions on the unit disk $\mathbb D\subset\mathbb C$ whose boundary values belong to the algebra $SAP(\partial \mathbb D) \subset L^\infty(\partial \mathbb D)$ of \textit{semi-almost periodic functions} on the unit circle $\partial\mathbb D$.  

A function $f \in L^\infty(\partial\mathbb D)$ is called semi-almost periodic if for any $s\in \partial \mathbb D$ and any $\varepsilon>0$ there exist functions $f_k:\partial \mathbb D \to\mathbb C$ ($k\in\{-1,1\}$) and arcs $\gamma_k$ with $s$ being their right (if $k=-1$) or left (if $k=1$) endpoint with respect to the counterclockwise orientation of $\partial\mathbb D$ such that the functions
$x \mapsto f_k\bigl(se^{i ke^x}\bigr)$, $-\infty<x<0$, $k\in\{-1,1\}$,
are restrictions of Bohr's almost periodic functions on $\mathbb R$ (see Definition \ref{def2.1} below)
and
\begin{equation*}
\sup_{z\in\gamma_k}|f(z)-f_k(z)|<\varepsilon , \quad k\in\{-1,1\}.
\end{equation*}

The graph of a real-valued semi-almost periodic function discontinuous at a single point has a form
\vspace*{-1mm}
\begin{center}
\includegraphics[scale=0.7]{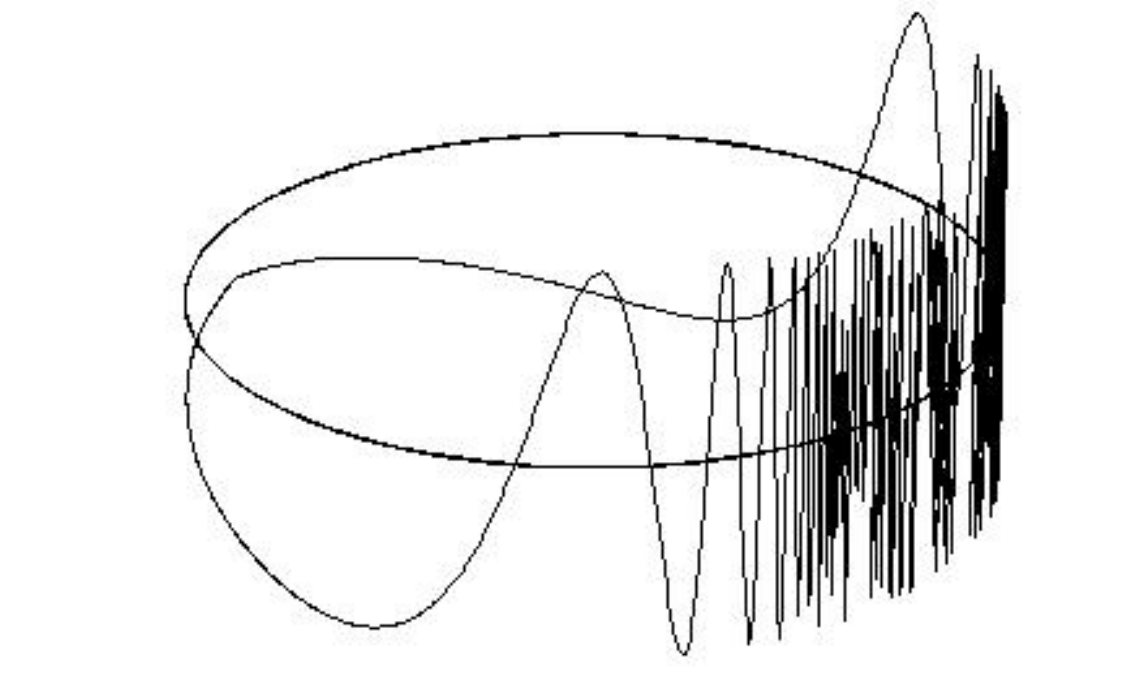}
\end{center}

Algebra $SAP(\partial \mathbb D)$ contains as a special case an algebra introduced by Sarason \cite{Sar}
in connection with some problems in the theory of Toeplitz operators. 
Our primary interest in holomorphic semi-almost periodic functions was motivated by the problem of description of the weakest possible boundary discontinuities of functions in $H^\infty(\mathbb D)$, the Hardy algebra of bounded holomorphic functions on $\mathbb D$. 
(Recall that a function $f \in H^\infty(\mathbb D)$ has radial limits almost everywhere on $\partial \mathbb D$,  the limit function $f|_{\partial \mathbb D} \in L^\infty(\partial \mathbb D)$, and $f$ can be recovered from $f|_{\partial \mathbb D}$ by means of the Cauchy integral formula.) 
In the general form this problem is as follows (see also \cite{BK1}): 

Given a continuous function $\Phi:\mathbb C \to \mathbb C$ to describe the minimal Banach subalgebra $H_{\Phi}^\infty(\mathbb D) \subset H^\infty(\mathbb D)$ containing all elements $f\in H^\infty(\mathbb D)^*$ such that $\Phi(f)|_{\partial \mathbb D}$ is piecewise Lipschitz having finitely many first-kind discontinuities. 

Here $H^\infty(\mathbb D)^*$ is the group of invertible elements of $H^\infty(\mathbb D)$.

Clearly, each $H_{\Phi}^\infty(\mathbb D)$ contains the {\em disk-algebra} $A(\mathbb D)$ (i.e., the algebra of holomorphic functions continuous up to the boundary). Moreover, if $\Phi(z)=z$, $\Real(z)$ or $\Imag(z)$, then the Lindel\"{o}f theorem, see, e.g., \cite{Gar}, implies that $H^\infty_{\Phi}(\mathbb D)=A(\mathbb D)$.
In contrast, if $\Phi$ is constant on a closed simple curve which does not encompass $0\in\mathbb C$, then $H^\infty_{\Phi}(\mathbb D)=H^\infty(\mathbb D)$. (This result is obtained by consequent applications of the Carath\'eodory conformal mapping theorem, the Mergelyan theorem and the Marshall theorem, see, e.g., \cite{Gar}.)
In \cite{BK} we studied the case of $\Phi(z)=|z|$ and showed that $H_{\Phi}^\infty(\mathbb D)$ coincides with the algebra of holomorphic semi-almost periodic functions $SAP(\partial \mathbb D) \cap H^\infty(\mathbb D)$.
In the present paper we continue the investigation started in \cite{BK}.
Despite the fact that our results concern the particular choice of $\Phi(z)=|z|$, the methods developed here and in \cite{BK} can be applied further to a more general class of functions $\Phi$.


Let $b\mathbb D$ be the maximal ideal space of the algebra $SAP(\partial \mathbb D) \cap H^\infty(\mathbb D)$, i.e., the set of all nonzero homomorphisms $SAP(\partial \mathbb D) \cap H^\infty(\mathbb D) \to \mathbb C$ equipped with the {\em Gelfand topology}. The disk $\mathbb D$ is naturally embedded into $b\mathbb D$. In \cite{BK} we proved that $\mathbb D$ is dense in $b\mathbb D$ (the so-called {\em corona theorem} for $SAP(\partial \mathbb D) \cap H^\infty(\mathbb D)$). We also described the topological structure of $b\mathbb D$.
In the present paper we refine and extend some of these results. In particular, we introduce Bohr-Fourier coefficients and spectra of functions from $SAP(\partial \mathbb D)$, describe \v{C}ech cohomology groups of $b\mathbb D$ and establish projective freeness of certain subalgebras of $SAP(\partial \mathbb D) \cap H^\infty(\mathbb D)$.
Recall that a commutative ring $R$ with identity is called {\em projective free} if every finitely generated projective $R$-module is free. Equivalently, $R$ is projective free iff every square idempotent matrix $F$ with entries in $R$ (i.e., such that $F^2=F$) is conjugate over $R$ to a matrix of the form
\begin{equation*}
\left(
\begin{array}{cc}
I_k&0\\
0&0
\end{array}
\right),
\end{equation*}
where $I_k$ stands for the $k \times k$ identity matrix. Every field $\mathbb F$ is trivially projective free. Quillen and Suslin proved  that if $R$ is projective free, then the rings of polynomials $R[x]$ and formal power series $R[[x]]$ over $R$ are projective free as well (see, e.g., \cite{QS}). 
Grauert proved that the ring ${\mathcal O}(\mathbb D^n)$ of holomorphic functions on the unit polydisk $\mathbb D^n$ is projective free \cite{Gra}. In turn, it was shown in \cite{BS} that the triviality of any complex vector bundle of finite rank over the connected maximal ideal space of a unital semi-simple commutative complex Banach algebra is sufficient for its projective freeness.
We employ this result to show that subalgebras of $SAP(\partial \mathbb D) \cap H^\infty(\mathbb D)$ whose elements have their spectra in non-negative or non-positive semi-groups are projective free.
Note that if a unital semi-simple commutative complex Banach algebra $A$ is projective free, then it is {\em Hermite}, i.e., every finitely generated stably free $A$-module is free. Equivalently, $A$ is Hermite iff any $k\times n$ matrix, $k<n$, with entries in $A$ having rank $k$ at each point of the maximal ideal space of $A$ can be extended to an invertible $n\times n$ matrix with entries in $A$, see \cite{Cohn}. (Here the values of elements of $A$ at points of the maximal ideal space are defined by means of the Gelfand transform.)

Finally, we prove that $SAP(\partial \mathbb D) \cap H^\infty(\mathbb D)$ has the \textit{approximation property}. (This result strengthen the approximation theorem of \cite{BK}.)
Recall that a Banach space $B$ is said to have the {\em approximation property} if for every compact set $K\subset B$ and every $\varepsilon>0$ there is an operator $T:B\to B$ of finite rank so that $$\|Tx-x\|_B <\varepsilon \quad \text{ for every } \quad x\in K.$$ 
(Throughout this paper all Banach spaces are assumed to be complex.)

Although it is strongly believed that the class of spaces with the approximation property includes practically all spaces which appear naturally in analysis, it is not known yet even for the space
$H^\infty(\mathbb D)$ (see, e.g., 
the paper of Bourgain and Reinov 
\cite{BR} for some results in this direction).
The first example of a space which fails to have the approximation property was constructed by Enflo \cite{E}. Since Enflo's work several other examples of such spaces were constructed, for the references see, e.g., \cite{L}.
Many problems of Banach space theory admit especially simple solutions if one of the spaces under consideration has the approximation property. One of such problems is the problem of determination whether given two Banach algebras $A \subset C(X)$, $B \subset C(Y)$ ($X$ and $Y$ are compact Hausdorff spaces) their {\em slice algebra}
\begin{equation*}
S(A,B):=\{f \in C(X \times Y): f(\cdot,y) \in A \text{ for all } y \in Y, f(x,\cdot) \in B \text{ for all }x \in X\}
\end{equation*}
coincides with $A \otimes B$, the closure in $C(X\times Y)$ of the symmetric tensor product of $A$ and $B$. For instance, this is true if either $A$ or $B$ have the approximation property. The latter is an immediate consequence of the following result of Grothendieck.

Let $A \subset C(X)$ be a closed subspace, $B$ be a Banach space and $A_B \subset C_B(X):=C(X,B)$ be the Banach space of all continuous $B$-valued functions $f$ such that $\varphi(f) \in A$ for any $\varphi \in B^*$. By $A\otimes B$ we denote completion of symmetric tensor product of $A$ and $B$ with respect to norm
\begin{equation}
\label{eq1}
\left\|\sum_{k=1}^m a_k\otimes b_k\right\|:=\sup_{x\in X} \left\|\sum_{k=1}^m a_k(x) b_k\right\|_B\ \ \ {\rm with}\ \ \ a_k\in A,\ b_k\in B.
\end{equation}

\vspace*{1mm}

\begin{theorem}[\cite{Gro}]
\label{equivthm}
The following statements are equivalent:

1) $A$ has the approximation property;

2) $A \otimes B=A_B$ for every Banach space $B$.
\end{theorem}

Our proof of the approximation property for $SAP(S) \cap H^\infty(\mathbb D)$ is based on Theorem \ref{equivthm} and on a variant of the approximation theorem in \cite{BK} for Banach-valued analogues of 
algebra $SAP(S) \cap H^\infty(\mathbb D)$.

The paper is organized as follows. Section \ref{sect1} is devoted to the algebra of semi-almost periodic functions $SAP(\partial \mathbb D)$. In Section \ref{sect2} we formulate our main results on the algebra of holomorphic semi-almost periodic functions $SAP(\partial \mathbb D) \cap H^\infty(\mathbb D)$. All proofs are presented in Section \ref{sect3}.

\vspace*{2mm}

The results of the present paper have been announced in \cite{BK1}.

\vspace*{2mm}

\noindent
{\em \bf Acknowledgment.} We are grateful to S.\,Favorov, S.\,Kislyakov and O.\,Reinov for useful discussions.

\section{Preliminaries on semi-almost periodic functions}

\label{sect1}

We first recall the definition of a Bohr almost periodic function on $\mathbb R$. In what follows, by $C_b(\mathbb R)$ we denote the algebra of bounded continuous functions on $\mathbb R$ endowed with $\sup$-norm.

\begin{definition}[see, e.g., \cite{Bes}]\label{def2.1}
A function $f\in C_b(\mathbb R)$ is said to be \textit{almost periodic} if the family of its translates $\{S_\tau f\}_{\tau \in \mathbb R}$\,, $S_\tau f(x):=f(x+\tau)$, $x \in\mathbb R$, is relatively compact in $C_b(\mathbb R)$. 
\end{definition}

The basic example of an almost periodic function is given by the formula
$$x \mapsto \sum_{l=1}^m c_l e^{i\lambda_l x}, \quad c_l \in \mathbb C, \quad \lambda_l \in \mathbb R.$$

Let $AP(\mathbb R)$ be the Banach algebra of almost periodic functions endowed with $\sup$-norm. 
The main characteristics of an almost periodic function $f  \in AP(\mathbb R)$ are its Bohr-Fourier coefficients $a_\lambda(f)$ and the spectrum $\speca(f)$ defined in terms of the mean value
\begin{equation}
\label{meanf}
M(f):=\lim_{T \to +\infty}\frac{1}{2T}\int_{-T}^T f(x)dx .
\end{equation}
Specifically,
\begin{equation}
\label{specf}
a_\lambda(f):=M(fe^{-i \lambda x}), \quad \lambda \in \mathbb R.
\end{equation}
Then $a_\lambda(f) \ne 0$ for at most countably many values of $\lambda$, see, e.g., \cite{Bes}. These values constitute the spectrum $\speca(f)$ of $f$. In particular, if $f=\sum_{l=1}^\infty c_le^{i\lambda_l x}$ ($c_l \ne 0$ and $\sum_{l=1}^\infty |c_l|<\infty$), then $\speca(f)=\{\lambda_1,\lambda_2,\dots\}$.

One of the main results of the theory of almost periodic functions states that each function $f \in AP(\mathbb R)$ can be uniformly approximated by functions of the form $\sum_{l=1}^m c_l e^{i\lambda_l x}$ with $\lambda_l \in \speca(f)$.

Let $\Gamma \subset \mathbb R$ be a unital additive semi-group (i.e., $0 \in \Gamma$). It follows easily from the cited approximation result that the space $AP_\Gamma(\mathbb R)$ of almost periodic functions with spectra in $\Gamma$ forms a unital Banach subalgebra of $AP(\mathbb R)$. We will use the following result.

\begin{theorem}
\label{prop3.1}
$AP_\Gamma(\mathbb R)$ has the approximation property.
\end{theorem}

Next, we recall the definition of a semi-almost periodic function on $\partial \mathbb D$ introduced in \cite{BK}.
In what follows, we consider $\partial\mathbb D$ with the counterclockwise orientation. For $s:=e^{it}$, $t \in [0,2\pi)$, let
\begin{equation}
\label{arcs}
\gamma_{s}^k(\delta):=\{se^{ikx} :\ 0 \leq x<\delta<2\pi\}, \quad k \in \{-1,1\},
\end{equation}
be two open arcs having $s$ as the right and the left endpoints (with respect to the orientation), respectively.

\begin{definition}[\cite{BK}]
\label{def1}
A function $f \in L^\infty(\partial\mathbb D)$ is called {\em semi-almost periodic} if for any $s \in \partial \mathbb D$, and any $\varepsilon>0$ there exist a number $\delta=\delta(s,\varepsilon) \in (0,\pi)$ and functions $f_k:\gamma_{s}^k(\delta)\to\mathbb C$, $k\in\{-1,1\}$, such that functions
\begin{equation*}
\tilde{f}_k(x):=f_k\bigl(se^{ik\delta e^x}\bigr), \quad -\infty<x<0, \quad k \in \{-1,1\},
\end{equation*}
are restrictions of some almost periodic functions from $AP(\mathbb R)$,
and
\begin{equation*}
\sup_{z \in \gamma_{s}^k(\delta)}|f(z)-f_k(z)|<\varepsilon ,\ \ \ k\in\{-1,1\}.
\end{equation*}
\end{definition}

By $SAP(\partial\mathbb D)$ we denote the Banach algebra of semi-almost periodic 
functions on $\partial\mathbb D$ endowed with $\sup$-norm. 
It is easy to see that the set of points of discontinuity of a function in $SAP(\partial \mathbb D)$ is at most countable.
For $S$ being a closed subset 
of $\partial\mathbb D$ we denote by $SAP(S)$ the Banach algebra of semi-almost periodic 
functions on $\partial\mathbb D$ that are continuous on $\partial\mathbb D\setminus S$. 
(Note that the Sarason algebra introduced in \cite{Sar} is isomorphic to $SAP(\{z_0\})$,  $z_0\in\partial\mathbb D$.)

\begin{example}[\cite{BK}]
\label{basicex}
A function $g$ defined on $\mathbb R\sqcup (\mathbb R+i\pi)$ is said to belong to the space $AP(\mathbb R\sqcup (\mathbb R+i\pi))$ if the functions $g(x)$ and $g(x+i\pi)$, $x\in \mathbb R$, belong to $AP(\mathbb R)$. The space $AP(\mathbb R\sqcup (\mathbb R+i\pi))$ is a function algebra (with respect to $\sup$-norm).  

Given $s \in \partial \mathbb D$ consider the map $\varphi_{s}:\partial \mathbb D \setminus \{-s\} \to \mathbb R$,
$\varphi_{s}(z):=\frac{2i(s-z)}{s+z}$, and define a linear isometric embedding $L_{s}:AP(\mathbb R\sqcup (\mathbb R+i\pi))\to L^\infty(\partial\mathbb D)$ by the formula
\begin{equation}
\label{formF}
(L_{s}g)(z):=
(g\circ\Log\circ \varphi_{s})(z),
\end{equation}
where $\Log(z):=\ln|z|+i{\rm Arg}(z)$, $z\in\mathbb C\setminus\mathbb R_{-}$, and ${\rm Arg}:\mathbb C\setminus\mathbb R_{-}\to (-\pi,\pi)$ stands for the principal branch of the multi-function ${\rm arg}$. Then the range of $L_{s}$ is a subspace of $SAP(\{-s,s\}).$
\end{example}

\begin{theorem}[\cite{BK}]
\label{propconeq}
For every $s\in\partial\mathbb D$ there exists a homomorphism of Banach algebras $E_{s}:SAP(\partial\mathbb D)\to AP(\mathbb R\sqcup (\mathbb R +i\pi))$
of norm $1$ such that for each $f \in SAP(\partial \mathbb D)$ the function $f-L_{s}(E_{s}f) \in SAP(\partial \mathbb D)$ is continuous and equal to $0$ at $s$. Moreover, any bounded linear operator $SAP(\partial\mathbb D)\to AP(\mathbb R\sqcup (\mathbb R +i\pi))$ satisfying this property coincides with $E_{s}$.
\end{theorem}

The functions $f_{-1,s}(x):=(E_{s}f)(x)$ and $f_{1,s}(x):=(E_{s}f)(x+i\pi)$, $x\in\mathbb R$,  are used to define the
left ($k=-1$) and the right ($k=1$) mean values $M_{s}^k(f)$ of a function $f \in SAP(\partial \mathbb D)$ over $s$ (cf. Remark \ref{equivway} below). Precisely, we put
$$M_{s}^k(f):=M(f_{k,s}).$$

Similarly, we define the left ($k=-1$) and the right ($k=1$) Bohr-Fourier coefficients and spectra of $f$ over $s$ by the formulas
\begin{equation*}
a_\lambda^k(f,s):=a_{\lambda}(f_{k,s})
\end{equation*}
and
\begin{equation*}
\speca_{s}^k(f):=\{\lambda \in \mathbb R:a_\lambda^k(f,s) \ne 0\}.
\end{equation*}
It follows immediate from the properties of the spectrum of an almost periodic function on $\mathbb R$ that $\speca_{s}^k(f)$ is at most countable.

Let $\Sigma:S \times \{-1,1\} \to 2^{\mathbb R}$ be a set-valued map which associates with each $s \in S$, $k \in \{-1,1\}$ a unital semi-group $\Sigma(s,k) \subset \mathbb R$.
By $SAP_{\,\Sigma}(S) \subset SAP(S)$ we denote the Banach algebra of semi-almost periodic functions $f$ with
$\speca_s^k(f) \subset \Sigma(s,k)$
for all $s \in S$, $k \in \{-1,1\}$. By the definition homomorphism $E_{s}$ of Theorem \ref{propconeq} sends each $f\in SAP_\Sigma(S)$ to the pair of functions $f_{-1,s}$, $f_{1,s}$ such that $f_{k,s}\in AP_{\Sigma(s,k)}$, $k\in\{-1,1\}$.

For a unital semi-group $\Gamma\subset\mathbb R$ by 
$b_\Gamma(\mathbb R)$ we denote the maximal ideal space of algebra $AP_\Gamma(\mathbb R)$. (E.g., for $\Gamma=\mathbb R$
the space $b\mathbb R:=b_{\mathbb R}(\mathbb R)$, commonly called the \textit{Bohr compactification} of $\mathbb R$, is a compact abelian topological group viewed as the inverse limit of compact finite-dimensional tori. The group $\mathbb R$ admits a canonical embedding into $b\mathbb R$ as a dense subgroup.)

Let $b^S_{\Sigma} (\partial \mathbb D)$ be the maximal ideal space of algebra $SAP_{\,\Sigma}(S)$ and
$r^S_\Sigma:b^S_{\Sigma} (\partial \mathbb D) \to \partial \mathbb D$ be the map transpose to the embedding $C(\partial \mathbb D) \hookrightarrow SAP_{\,\Sigma}(S)$.
The proof of the next statement is analogous to the proof of Theorem 1.7 in \cite{BK}.

\begin{theorem}
\label{sapthm} 
\begin{itemize}
\item[(1)]
The map transpose to the restriction of homomorphism $E_{s}$ to $SAP_\Sigma(S)$ determines
an embedding $h_\Sigma^s: b_{\Sigma(s,-1)}(\mathbb R) \sqcup b_{\Sigma(s,1)}(\mathbb R)\hookrightarrow b^S_{\Sigma} (\partial \mathbb D)$ whose image coincides with $(r^S_\Sigma)^{-1}(s)$. 
\item[(2)]
The restriction $r^S_\Sigma:b^S_{\Sigma} (\partial \mathbb D) \setminus (r^S_{\Sigma})^{-1}(S) \to \partial \mathbb D \setminus S$ is a homeomorphism.
\end{itemize}
\begin{center}
\vspace*{-2mm}
\includegraphics[scale=0.6]{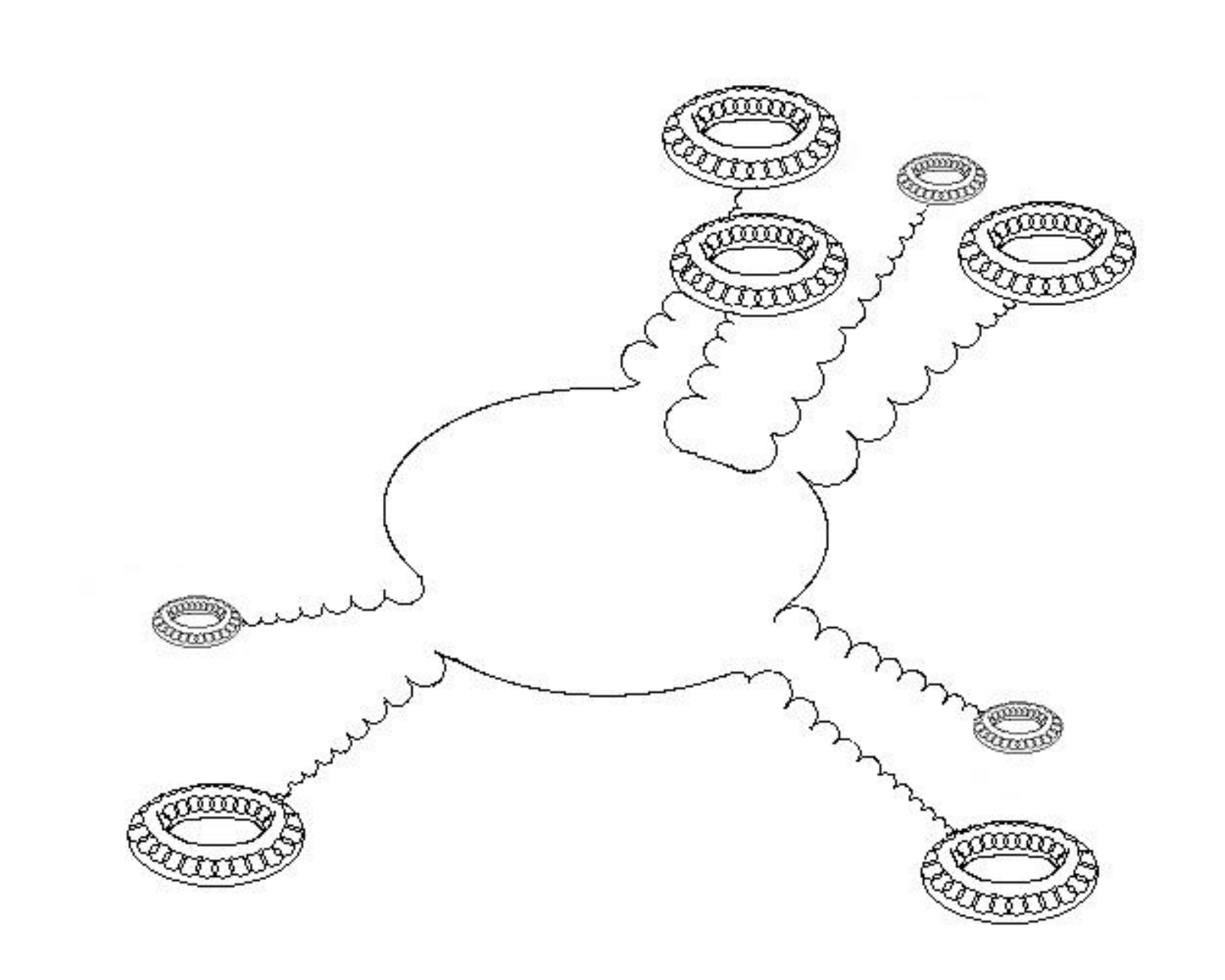}
\end{center}
\vspace*{-4mm}
\textit{(For an $m$ point set $S$ and each $\Sigma(s,k)$, $s\in S, \, k\in\{-1,1\}$, being a group, the maximal ideal space $b^S_\Sigma(\partial \mathbb D)$ is the union of $\partial \mathbb D \setminus S$ and $2m$ Bohr compactifications $b_{\Sigma(s,k)}(\mathbb R)$ that can be viewed as (finite or infinite dimensional) tori.)}
\end{theorem}

\begin{remark}
\label{equivway}
There is an equivalent way to define the mean value of a semi-almost periodic function.
Specifically, it is easily seen that for a semi-almost periodic function $f \in SAP(\partial \mathbb D)$, $k \in \{-1,1\}$, and a point $s \in S$ the left ($k=-1$) and the right ($k=1$) mean values of $f$ over $s$ are given by the formulas
\begin{equation*}
M_{s}^k(f):=\lim_{n \to \infty} \frac{1}{b_n-a_n}\int_{a_n}^{b_n}f\left (se^{ik e^{t}}\right)dt,
\end{equation*}
where $\{a_n\}$, $\{b_n\}$ are arbitrary sequences converging to $-\infty$ such that $\lim_{n\to\infty}(b_n-a_n)=+\infty$.

The Bohr-Fourier coefficients of $f$ over $s$ can be then defined by the formulas
\begin{equation*}
a_\lambda^k(f,s):=M_s^k(fe^{-i\lambda \log_s^k}),
\end{equation*}
where 
$$\log_s^k(se^{ikx}):=\ln x, \quad 0<x<2\pi, \quad k\in\{-1,1\}.$$
\end{remark}

The next result encompasses the basic properties of the mean value and the spectrum of a semi-almost periodic function.

\begin{theorem}
\label{totalspecprop} 
\begin{itemize}
\item[(1)] For each $s\in S$, $k\in\{-1,1\}$ the mean value $M_s^k$ is a complex continuous linear functional on $SAP(\partial \mathbb D)$ of norm $1$.
\item[(2)] A function $f\in SAP_{\,\Sigma}(S)$ if and only if for each $s \in S$ and $k \in \{-1,1\}$ the almost periodic functions $\tilde{f}_k$ in Definition \ref{def1} can be chosen from $AP_{\,\Sigma(s,k)}(\mathbb R)$.
\item[(3)] The ``total spectrum`` $\bigcup_{s \in S,\,k=\pm 1}\speca_s^k(f)$ of a function $f \in SAP(S)$ is at most countable.
\end{itemize}
\end{theorem}

\section{Holomorphic semi-almost periodic functions: main results}

\label{sect2}

{\bf 3.1.} Let $C_b(T)$ denote the complex Banach space of bounded continuous functions on the strip $T:=\{z\in\mathbb C\ :\ \Imag(z)\in [0,\pi]\}$ endowed with $\sup$-norm.

\begin{definition}[see, e.g., \cite{Bes}]
\label{def4}
A function $f\in C_b(T)$ is called \textit{holomorphic almost periodic} if it is holomorphic in the interior of $T$
and the family of its translates $\{S_x f\}_{x \in \mathbb R}$, $S_xf(z):=f(z+x)$, $z\in T$, is relatively compact in $C_b(T)$. 
\end{definition}

We denote by $APH(T)$ the Banach algebra of holomorphic almost periodic functions endowed with $\sup$-norm. Any function in $APH(T)$ is uniformly continuous on $T$.
The mean value of a function $f \in APH(T)$ is defined by the formula
\begin{equation}
\label{meanf2}
M(f):=\lim_{T \to +\infty}\frac{1}{2T}\int_{-T}^T f(x+iy)dx \in \mathbb C
\end{equation}
($M(f)$ does not depend on $y$, see, e.g., \cite{Bes}). Further, the Bohr-Fourier coefficients of $f$ are defined by
\begin{equation}
\label{specf2}
a_\lambda(f):=M(fe^{-i \lambda z}), \quad \lambda \in \mathbb R.
\end{equation}
Then $a_\lambda(f) \ne 0$ for at most countably many values of $\lambda$, these values form the spectrum $\speca(f)$ of $f$. For instance, if $f=\sum_{l=1}^\infty c_le^{i\lambda_l z}$ ($c_l \ne 0$, $\sum_{l=1}^\infty |c_l|<\infty$), then $\speca(f)=\{\lambda_1,\lambda_2,\dots\}$. 
Similarly to the case of functions from $AP(\mathbb R)$ each $f \in APH(T)$ can be uniformly approximated by functions of the form $\sum_{l=1}^m c_l e^{i\lambda_l z}$ with $\lambda_l \in \speca(f)$.

Let $\Gamma \subset \mathbb R$ be a unital additive semi-group. The space $APH_\Gamma(T)$ of holomorphic almost periodic functions with spectra in $\Gamma$ forms a unital Banach algebra. Analogously to Theorem \ref{prop3.1} one has

\begin{theorem}
\label{prop3.2}
$APH_\Gamma(T)$ has the approximation property.
\end{theorem}

The functions in $SAP_\Sigma(S) \cap H^\infty(\mathbb D)$ are called {\em holomorphic semi-almost periodic}.

\begin{example}[\cite{BK}]
\label{ex2}
For $s \in \partial \mathbb D$ consider the map
\begin{equation*}
\varphi_s:\bar{\mathbb D} \setminus \{-s\} \to \bar{\mathbb H}_+, \quad \varphi_s(z):=\frac{2i(s-z)}{s+z}.
\end{equation*}
Here $\mathbb H^{+}$ is the upper half-plane. Then $\varphi_s$ maps $\mathbb D$ conformally onto $\mathbb H_+$ and $\partial\mathbb D\setminus\{-s\}$ diffeomorphically onto $\mathbb R$ (the boundary of $\mathbb H_{+}$) so that $\varphi_{s}(s)=0$.

Let $T_{0}$ be the interior of the strip $T$. Consider the conformal map ${\rm Log}:\mathbb H_{+}\to T_{0}$,
$z\mapsto {\rm Log}(z):=\ln|z|+i{\rm Arg}(z)$, where ${\rm Arg}:\mathbb C\setminus\mathbb R_{-}\to (-\pi,\pi)$ is the principal branch of the multi-function ${\rm arg}$. The function ${\rm Log}$ is extended to a homeomorphism of $\overline{\mathbb H}_{+}\setminus\{0\}$ onto $T$. 
Let $g\in APH(T)$. Then the function
\begin{equation*}
(L_{s}g)(z):=(g \circ \Log\circ \varphi_s)(z), \quad z \in \mathbb D,
\end{equation*}
belongs to $SAP(\{-s,s\}) \cap H^\infty(\mathbb D)$.
\end{example}

\begin{proposition}
\label{impthm}
Suppose that $f \in SAP(S)\cap H^\infty(\mathbb D)$. Then $\speca_s^{-1}(f)= \speca_s^1(f)=:\speca_s(f)$ and, moreover, 
\begin{equation}
\label{speccond3}
a^{-1}_\lambda(f,s)=e^{\lambda \pi}a^{1}_\lambda(f,s) \quad \text{for each} \quad  \lambda\in \speca_s(f).
\end{equation}
\end{proposition}
\noindent (Recall that the choice of the upper indices $\pm 1$ is determined by the orientation of $\partial \mathbb D$.)

Proposition \ref{impthm} and Lindel\"{o}f's theorem (see, e.g., \cite{Gar}) imply that 
$SAP_\Sigma(S) \cap H^\infty(\mathbb D)=SAP_{\Sigma'}(S') \cap H^\infty(\mathbb D),$
where 
$$S':=\{s \in S: \Sigma(s,-1) \cap \Sigma(s,1) \ne \{0\}\}$$
and
$$\Sigma'(s,k):=\Sigma(s,-1) \cap \Sigma(s,1) \quad \text{ for }\quad k=-1,1, \quad s \in S'.$$ 

\begin{quote}
In what follows we assume that
$\Sigma(s,-1)=\Sigma(s,1)=:\Sigma(s)$ and each $\Sigma(s)$, $s \in S$, is non-trivial.
\end{quote}

\begin{example}
\label{ex5}
If $g(z):=e^{\frac{i \lambda}{\pi}z}$, $z\in T$, then $L_{s}g=e^{\lambda h}$, where $h$ is a holomorphic functions whose real part $\Real(h)$ is the characteristic function of the closed arc going in the counterclockwise direction from the initial point at $s$ to the endpoint at $-s$, and such that $h(0)=\frac 12+\frac{i\ln 2}{\pi}$. Thus 
$$\speca_{s}(e^{\lambda h})=\{\lambda/\pi\}.$$
Indeed, in this case the restriction of $g\circ \Log$ to $\mathbb R$ is equal to $x \mapsto e^{\lambda(\chi_{\mathbb R_+}(x)+\frac{i\ln|x|}{\pi})}$,
where $\chi_{\mathbb R_+}$ is the characteristic function of $\mathbb R_+$.
In turn, the restriction of the pullback $e^{\lambda h} \circ \varphi_{s}^{-1}$ to $\mathbb R$ coincides with $e^{\lambda(\chi_{\mathbb R_+}+\frac{i \ln|x|}{\pi})}$ as well. This implies the required result.
\end{example}

\noindent {\bf 3.2.}
The main result of this section is

\begin{theorem}
\label{semithm}
$SAP_\Sigma(S) \cap H^\infty(\mathbb D)$ has the approximation property.
\end{theorem}

Our proof of Theorem \ref{semithm} is based on the equivalence established in Theorem \ref{equivthm} and on an approximation result for  Banach-valued analogues of 
algebra $SAP(S) \cap H^\infty(\mathbb D)$ formulated below. Specifically, for a Banach space $B$ we define
\begin{equation*}
SAP^B_\Sigma(S):=SAP_\Sigma(S) \otimes B.
\end{equation*}
Using the Poisson integral formula we can extend each function from $SAP^B_\Sigma(S)$ to a bounded $B$-valued harmonic function on $\mathbb D$ having the same $\sup$-norm. We identify $SAP^B_\Sigma(S)$ with its harmonic extension.
Let $H^\infty_B(\mathbb D)$ be the Banach space of bounded $B$-valued 
holomorphic functions on $\mathbb D$  equipped with $\sup$-norm.
By $(SAP_\Sigma(S)\cap H^\infty(\mathbb D))_B$ we denote the Banach space of all continuous $B$-valued functions $f$ on the maximal ideal space $b^S(\mathbb D)$ of algebra $SAP(S)\cap H^\infty(\mathbb D)$ such that 
$\varphi(f)\in SAP_\Sigma(S) \cap H^\infty(\mathbb D)$ for any $\varphi \in B^*$. In what follows we naturally identify
$\mathbb D$ with a subset of $b^S(\mathbb D)$.

\begin{proposition}
\label{strprop}
Let $f\in (SAP_\Sigma(S)\cap H^\infty(\mathbb D))_B$. Then $f|_{\mathbb D}\in SAP^B_\Sigma(S) \cap H^\infty_B(\mathbb D)$. 
\end{proposition}

Let $A_\Sigma^S$ be the closed subalgebra of $H^\infty(\mathbb D)$ generated by the disk-algebra $A(\mathbb D)$ and the functions of the form $ge^{\lambda h}$,
where $\Real(h)|_{\partial\mathbb D}$ is the characteristic function of the closed arc going in the counterclockwise direction from the initial point at $s$ to the endpoint at $-s$ such that $s\in S$, $\frac{\lambda}{\pi} \in \Sigma(s)$ and $g(z):=z+s$, $z\in\mathbb D$ (in particular, $ge^{\lambda h}$ has discontinuity at $s$ only).

The next result combined with Proposition \ref{strprop} and Theorem \ref{equivthm} implies Theorem \ref{semithm}.

\begin{theorem}
\label{Bapprox}
$SAP^B_\Sigma(S) \cap H^\infty_B(\mathbb D)=A_\Sigma^S \otimes B$.
\end{theorem}
As a corollary we obtain
\begin{corollary}
\label{thm4}
\label{oldapprox}
$SAP_\Sigma (S) \cap H^\infty(\mathbb D)=A_\Sigma^S$.
\end{corollary}


This immediately implies the following result.
\begin{theorem}\label{appr}
$SAP_\Sigma(S) \cap H^\infty(\mathbb D)$ is generated by algebras $SAP_{\Sigma|_F}(F) \cap H^\infty(\mathbb D)$ for all possible finite subsets $F$ of $S$.
\end{theorem}

\vspace*{2mm}


\noindent\textbf{3.3.} 
The algebras $SAP_\Sigma(S) \cap H^\infty(\mathbb D)$ are preserved under the action of the group $\Aut(\mathbb D)$ of biholomorphic automorphisms $\mathbb D \to \mathbb D$. More precisely, each $\kappa\in \Aut(\mathbb D)$ is extended to a diffeomorphism $\bar{\mathbb D}\to\bar{\mathbb D}$ (denoted by the same symbol). 
We denote by $\kappa^*:H^\infty(\mathbb D) \to H^\infty(\mathbb D)$ the pullback by $\kappa$, and put $\kappa_* S:=\kappa(S)$, $(\kappa_*\Sigma)(s,\cdot):=\Sigma(\kappa(s),\cdot).$ 

\begin{theorem}
\label{isothm}
$\kappa^*$ maps $SAP_{\kappa_*\Sigma}(\kappa_* S) \cap H^\infty(\mathbb D)$ isometrically and isomorphically onto  $SAP_\Sigma(S) \cap H^\infty(\mathbb D)$.
\end{theorem}

We conclude this section with a result on the tangential behavior of functions from $SAP(\partial \mathbb D) \cap H^\infty(\mathbb D)$.

\begin{theorem}
\label{convthm}
Let $\{z_n\}_{n\in\mathbb N} \subset \mathbb D$ and $\{s_n\}_{n\in\mathbb N} \subset \partial \mathbb D$ converge to a point $s_0 \in \partial \mathbb D$. Assume that
\begin{equation}
\label{convcond}
\lim_{n\to\infty}\frac{|z_n-s_n|}{|s_0-s_n|}=0.
\end{equation}
Then for every $f \in SAP(\partial \mathbb D) \cap H^\infty(\mathbb D)$ the limits 
$\lim_{n\to\infty} f(z_n)$ and $\lim_{n\to\infty} f(s_n)$ do not exist or exist simultaneously and in the latter case they are equal. 
\end{theorem}

\begin{remark}
This result implies
that the extension (by means of the Gelfand transform) of each $f \in SAP(\partial \mathbb D) \cap H^\infty(\mathbb D)$ to the maximal ideal space of $H^\infty(\mathbb D)$ is constant on a nontrivial Gleason part containing a limit point of a net in $\mathbb D$ converging tangentially to $\partial\mathbb D$. In turn, one can easily show that if $s\in S$ and the minimal subgroup of $\mathbb R$ containing $\Sigma(s)$ is not isomorphic to $\mathbb Z$, then $SAP_{\Sigma}(S)\cap H^\infty(\mathbb D)$ separates points of each nontrivial Gleason part containing a limit point of a net in $\mathbb D$ converging non-tangentially to $s$ (we refer to \cite{Gar} for the corresponding definitions).
\end{remark}

In the next two sections we formulate some topological results about the maximal ideal spaces of algebras $SAP_\Sigma(S) \cap H^\infty(\mathbb D)$.

\vspace*{2mm}

\noindent\textbf{3.4.~}Let $b^S_\Sigma(\mathbb D)$ denote the maximal ideal space of $SAP_\Sigma(S) \cap H^\infty(\mathbb D)$.
The inclusion $$SAP_{\Sigma|_{F_1}}(F_1) \cap H^\infty(\mathbb D) \subset SAP_{\Sigma|_{F_2}}(F_2) \cap H^\infty(\mathbb D) \quad \text{ if } \quad F_1 \subset F_2$$ determines a continuous map of maximal ideal spaces $$\omega_{F_1}^{F_2}: b^{F_2}_{\Sigma|_{F_2}}(\mathbb D) \to b^{F_1}_{\Sigma|_{F_1}}(\mathbb D).$$ The family $\{b^F_{\Sigma|_F}(\mathbb D)\,;\,\omega\}_{F\subset S\,;\, \#F<\infty}$  forms the inverse limiting system. From Theorem \ref{appr} we obtain

\begin{theorem}
\label{invlimcor}
$b^S_\Sigma(\mathbb D)$ is the inverse limit of $\{b^F_{\Sigma|_F}(\mathbb D)\,;\,\omega\}_{F\subset S\,;\, \#F<\infty}$.
\end{theorem}

Let
\begin{equation}\label{dualmap}
a^S_\Sigma:b^S_\Sigma(\mathbb D) \to \bar{\mathbb D}
\end{equation}
be the continuous surjective map transpose to the embedding $A(\mathbb D) \hookrightarrow SAP_\Sigma(S) \cap H^\infty(\mathbb D)$. (Recall that the maximal ideal space of the disk-algebra $A(\mathbb D)$ is homeomorphic to $\bar{\mathbb D}$.) By $b_\Gamma(T)$ we denote the maximal ideal space of algebra $APH_\Gamma(T)$ and by $\iota_\Gamma:T\to b_\Gamma(T)$ the continuous map determined by evaluations at points of $T$.
(Observe that $\iota_\Gamma$ is not necessarily an embedding.)

The proof of the next statement can be obtained by following closely the arguments in the proof of Theorem 1.14 in \cite{BK}. In its formulation we
assume that the corresponding algebras are defined on their maximal ideal spaces by means of the Gelfand transforms.
\begin{theorem} 
\label{maxidthm}
\begin{itemize}
\item[(1)] 
For each $s \in S$ there exists an embedding $i_\Sigma^s: b_{\Sigma(s)}(T)\hookrightarrow b_\Sigma^S(\mathbb D)$ whose image is
$(a^S_\Sigma)^{-1}(s)$ such that the pullback $(i_\Sigma^s)^*$ maps $SAP_\Sigma(S) \cap H^\infty(\mathbb D)$ surjectively onto  $APH_{\Sigma(s)}(T)$. Moreover, the composition of the restriction map to $\mathbb R\sqcup (\mathbb R+i\pi)$ and $(i_\Sigma^s\circ\iota_{\Sigma(s)})^*$ coincides with
the restriction of homomorphism $E_{s}$ to $SAP_\Sigma(S) \cap H^\infty(\mathbb D)$ (see Theorem \ref{propconeq}).
\item[(2)]
The restriction $a^S_\Sigma: b^S_\Sigma(\mathbb D) \setminus (a^S_\Sigma)^{-1}(S) \to \bar{\mathbb D} \setminus S$ is a homeomorphism.
\end{itemize}
\end{theorem}

Since $SAP_\Sigma(S)\cap H^\infty(\mathbb D)$ separates the points on $\mathbb D$, the evaluation at points of $\mathbb D$ determines a natural embedding $\iota:\mathbb D \hookrightarrow b^S_\Sigma(\mathbb D)$. 

One has the following commutative diagram of maximal ideal spaces considered in the present paper, where the `dashed' arrows stand for embeddings in the case $\Sigma(s,-1)=\Sigma(s,1)$ are (non-trivial) groups for all $s \in S$, and for continuous maps otherwise.
\begin{center}
\vspace*{2mm}
\includegraphics[scale=0.65]{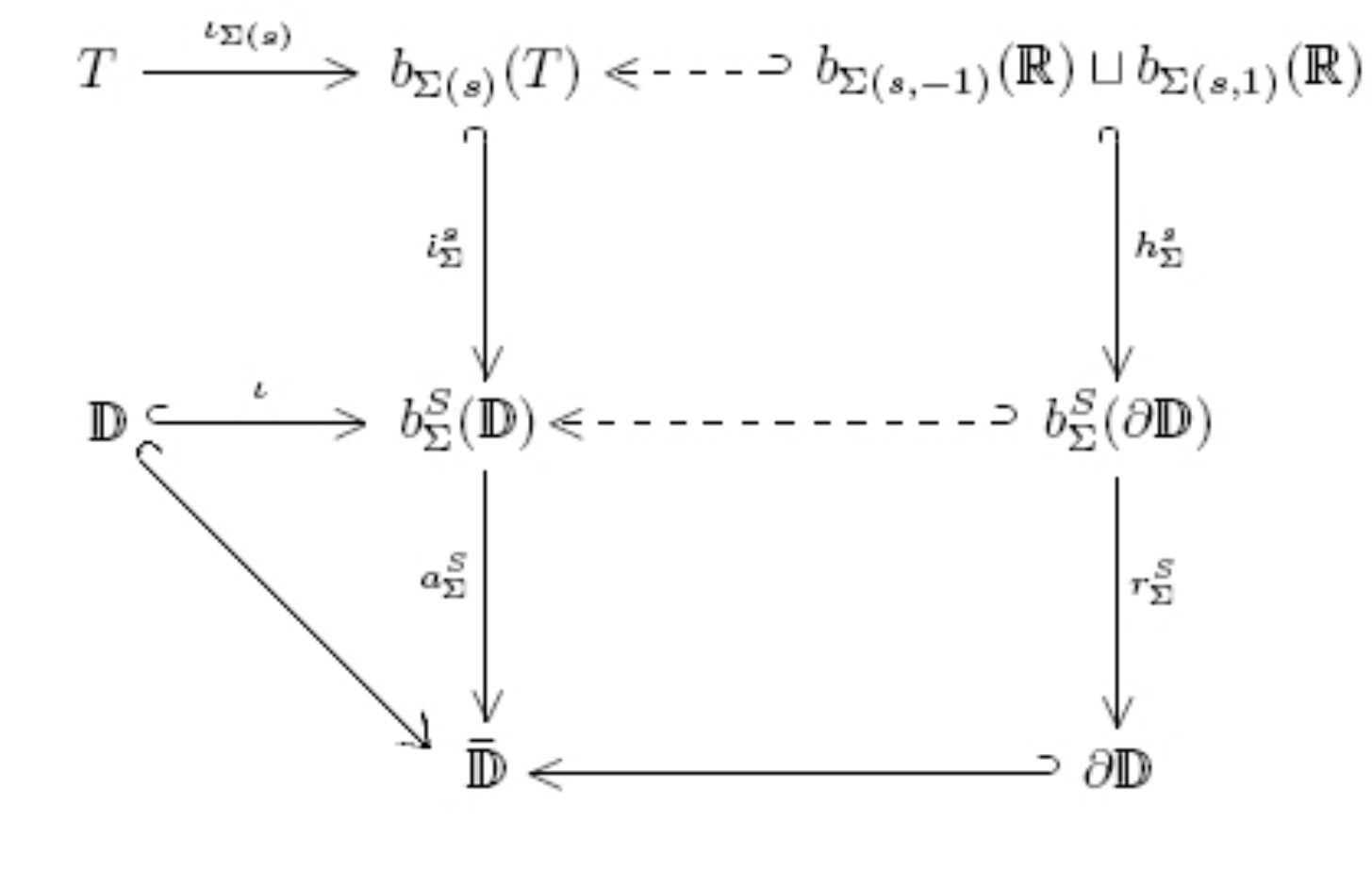}
\end{center}
\begin{theorem}[Corona Theorem]
\label{coronathm}
$\iota(\mathbb D)$ is dense in $b^S_\Sigma(\mathbb D)$ iff each $\Sigma(s)$, $s\in S$, is a group.
\end{theorem}

Recall that the corona theorem is equivalent to the following statement: {\em for any collection of functions $f_1, \dots, f_m \in SAP_\Sigma(S) \cap H^\infty(\mathbb D)$ such that $$\max_{1 \leq k \leq m}|f_k(z)| \geq \delta>0 \quad \text{ for all } z \in \mathbb D$$ there exist functions $g_1,\dots,g_m \in SAP_\Sigma(S) \cap H^\infty(\mathbb D)$ such that}
\begin{equation*}
f_1g_1+\dots+f_mg_m =1 \quad \text{ on } \mathbb D.
\end{equation*}

Our next result shows that $b_\Sigma^S(\mathbb D)$, $S\neq\emptyset$, is not arcwise connected.

\begin{theorem}
\label{holmapsthm}
Assume that $F: [0,1]\to b_{\Sigma}^S(\mathbb D)$ is continuous. Then either $F([0,1]) \subset \bar{\mathbb D} \setminus S$ or there exists $s \in S$ such that $F([0,1]) \subset (a_{\Sigma}^S)^{-1}(s)$.
\end{theorem}
\begin{remark}
From Theorem \ref{holmapsthm} one obtains straightforwardly a similar statement with $[0,1]$ replaced by an arcwise connected topological space.
\end{remark}

\noindent\textbf{3.5.~}Let $K^S_\Sigma$ be the \v{S}ilov boundary of algebra $SAP_\Sigma(S) \cap H^\infty(\mathbb D)$, that is, the minimal closed subset of $b^S_\Sigma(\mathbb D)$ such that for every $f \in SAP_\Sigma(S) \cap H^\infty(\mathbb D)$
\begin{equation*}
\sup_{z \in \mathbb D}|f(z)|=\max_{\varphi \in K^S_\Sigma}|f(\varphi)|,
\end{equation*}
where $f$ is assumed to be extended to $b^S_\Sigma(\mathbb D)$ by means of the Gelfand transform.
For a non-trivial semigroup $\Gamma\subset\mathbb R$ by $\cl_{\Gamma}(\mathbb R+i\pi)$ and $\cl_{\Gamma}(\mathbb R)$ we denote closures of $\iota_\Gamma(\mathbb R+i\pi)$ and $\iota_\Gamma(\mathbb R)$ in $b_{\Gamma}(T)$\, (the maximal ideal space of 
$APH_{\Gamma}(T)$). One can easily show that these closures are homeomorphic to $b_{\widehat\Gamma}(\mathbb R)$, where $\widehat\Gamma$ is the minimal subgroup of $\mathbb R$ containing $\Gamma$.

We retain notation of Theorem \ref{maxidthm}.

\begin{theorem}

\label{shilovthm}
$$
K^S_\Sigma=\left(\bigcup_{s\in S} i_\Sigma^s\bigl(\cl_{\Sigma(s)}(\mathbb R) \cup \cl_{\Sigma(s)}(\mathbb R+i\pi)\bigr)\right) \cup \partial \mathbb D \setminus S .
$$
\end{theorem}

\begin{remark}
If each $\Sigma(s)$, $s\in S$, is a group, then the \v{S}ilov boundary $K_\Sigma^S$ is naturally homeomorphic to the maximal ideal space $b_\Sigma^S(\partial\mathbb D)$ of algebra $SAP_\Sigma(S)$, cf. Theorem \ref{sapthm}.
\end{remark}

Next, we formulate a result on the \v{C}ech cohomology groups of $b^S_\Sigma(\mathbb D)$.

\begin{theorem} 
\label{cechthm} 
\begin{itemize}
\item[(1)] The \v{C}ech cohomology groups
$$H^k\bigl(b^S_\Sigma(\mathbb D),\mathbb Z\bigr)\cong\bigoplus_{s \in S}H^k\bigl(b_{\Sigma(s)}(T),\mathbb Z\bigr), \quad k \geq 1.$$
\item[(2)] Suppose that each $\Sigma(s)$ is a subset of $\mathbb R_+$ or $\mathbb R_{-}$. Then
$H^k\bigl(b^S_\Sigma(\mathbb D),\mathbb Z\bigr)=0$, $k \geq 1$, and $SAP_\Sigma(S) \cap H^\infty(\mathbb D)$ is projective free.
\end{itemize}
\end{theorem}

If $G$ is a compact connected abelian topological group and $\hat{G}$ is its dual, then $H^k(G,\mathbb Z) \cong \wedge_{\mathbb Z}^k \hat{G}$, $k \geq 1$ (see, e.g., \cite{HM}).  
Using Pontryagin duality one obtains

\begin{corollary}
\label{finalcor}
Assume that each $\Sigma(s)$ is a group. Then 
$$H^k(b_{\Sigma}^S(\mathbb D),\mathbb Z)\cong\bigoplus_{s \in S}\wedge_{\mathbb Z}^k \Sigma(s), \quad k \geq 1.$$
In particular, if for a fixed $n\in\mathbb N$ each $\Sigma(s)$ is isomorphic to a subgroup of $\mathbb Q^{n}$, then $H^k(b_{\Sigma}^S(\mathbb D),\mathbb Z)=0$ for all $k \geq n+1$.
\end{corollary}

Finally, we describe the set of connected components of the group of invertible matrices with entries in $SAP_\Sigma(S) \cap H^\infty(\mathbb D)$.

Let $GL_n(A)$ denote the group of invertible $n \times n$ matrices with entries in a unital Banach algebra $A$. By $[GL_n(A)]$ we denote the group of  connected components of $GL_n(A)$, i.e., the quotient of $GL_n(A)$ by the connected component containing the unit $I_n\in GL_n(A)$ (this is clearly a  normal subgroup of $GL_n(A)$).

We set $$G^n_{\Sigma(s)}(T):=GL_n\bigl(APH_{\Sigma(s)}(T)\bigr)\quad  \text{and}\quad G^n_\Sigma(S):=GL_n\bigl(SAP_\Sigma(S) \cap H^\infty(\mathbb D)\bigr).$$

Let $bT(S):=\sqcup_{s\in S}\, b_{\Sigma(s)}(T)$. According to Theorem \ref{maxidthm} there exists a natural embedding
$i_\Sigma^s: b_{\Sigma(s)}(T)\hookrightarrow b_\Sigma^S(\mathbb D)$ whose image is $(a^S_\Sigma)^{-1}(s)$. Then the map
$$
I: bT(S)\to (a_\Sigma^S)^{-1}(S),\qquad I(\xi):=i_\Sigma^s(\xi)\quad \text{for}\quad \xi \in b_{\Sigma(s)}(T),
$$
is a bijection.

\begin{theorem}
\label{concompthm}
The map transpose to the composition $bT(S)\stackrel{I}{\longrightarrow}(a^S_\Sigma)^{-1}(S)\hookrightarrow b_\Sigma^S(\mathbb D)$ induces an isomorphism
$$[G^n_\Sigma(S)] \cong \bigoplus_{s \in S} [G^n_{\Sigma(s)}(T)].$$
In particular, if each $\Sigma(s)$, $s\in S$, is a subset of $\mathbb R_+$ or $\mathbb R_{-}$, then $G_\Sigma^n(S)$ is connected.
\end{theorem}

\begin{remark}
According to a result of Arens \cite{A}, $[G^n_{\Sigma(s)}(T)]$, $s\in S$, can be identified with the group $[b_{\Sigma(s)}(T), GL_n(\mathbb C)]$ of homotopy classes of continuous maps $b_{\Sigma(s)}(T)\to GL_n(\mathbb C)$. Moreover, if $\Sigma(s)$ is a group, then $b_{\Sigma(s)}(T)$ is homotopically equivalent to $b_{\Sigma(s)}(\mathbb R)$, the maximal ideal space of algebra $AP_{\Sigma(s)}(\mathbb R)$ (see the proof of Corollary \ref{finalcor}). In this case $b_{\Sigma(s)}(\mathbb R)$ is the inverse limit of a family of finite-dimensional tori. Then $[b_{\Sigma(s)}(T), GL_n(\mathbb C)]$ is isomorphic to the direct limit of torus homotopy groups with values in $GL_n(\mathbb C)$ corresponding to this family. As follows from the classical results of Fox \cite{F}, the latter can be expressed as a direct sum of certain homotopy groups of the unitary group $U_n\subset GL_n(\mathbb C)$.
\end{remark}

\section{Proofs}

\label{sect3}

\subsection{Proofs of Theorems \ref{prop3.1} and \ref{prop3.2}}

We will prove Theorem \ref{prop3.2} only (the proof of Theorem \ref{prop3.1} is similar). We refer to the book of Besicovich \cite{Bes} for the corresponding definitions and facts from the theory of almost periodic functions.

\begin{proof}
Let $K \subset APH_\Gamma(T)$ be compact. Given $\varepsilon>0$ consider an $\frac{\varepsilon}{3}$-net $\{f_1,\dots, f_l\}\subset K$. Let 
\begin{equation*}
\mathcal K(t):=\sum_{|\nu_1| \leq n_1, \dots, |\nu_r| \leq n_r} \left(1-\frac{\nu_1}{n_1} \right)\dots \left(1-\frac{\nu_r}{n_r} \right)e^{-i \left(\frac{\nu_1}{n_1}\beta_1+\dots+\frac{\nu_r}{n_r}\beta_r \right)t}
\end{equation*}
be a Bochner-Fejer kernel such that for all $1\le k\le l$
\begin{equation}\label{eq3.1}
\sup_{z\in\Sigma}|f_k(z)-M_t\{f_k(z+t)\mathcal K(t)\}|\le\frac{\varepsilon}{3}.
\end{equation}
Here $\beta_1,\dots,\beta_r$ are linearly independent over $\mathbb Q$ and belong to the union of spectra of functions $f_1,\dots, f_l$,\, $\nu_1,\dots, \nu_r\in\mathbb Z$, $n_1,\dots, n_r\in\mathbb N$, and
$$
M_t\{f_k(z+t)\mathcal K(t)\}:=\lim_{T\to\infty}\frac{1}{2T}\int_{-T}^{T}f_k(z+t)\mathcal K(t)\,dt
$$
are the corresponding holomorphic Bochner-Fejer polynomials belonging to $APH_\Gamma(T)$ as well (clearly, the spectrum of the function $z \mapsto M_t\{f_k(z+t)\mathcal K(t)\}$ is contained in $\speca(f_k)$).

We define a linear operator $T: APH_\Gamma(T)\to APH_\Gamma(T)$ from the definition of the approximation property  by the formula
\begin{equation}\label{eq3.2}
(Tf)(z):=M_t\{f(z+t)\mathcal K(t)\},\ \ \ f\in APH_\Gamma(T).
\end{equation}
Then $T$ is a bounded linear projection onto a finite-dimensional subspace of $APH_\Gamma(T)$ generated by functions $e^{i \left(\frac{\nu_1}{n_1}\beta_1+\dots+\frac{\nu_r}{n_r}\beta_r \right)z}$, $|\nu_1| \leq n_1, \dots, |\nu_r| \leq n_r$. Moreover, since $\mathcal K(t)\ge 0$ for all $t\in\mathbb R$ and
$M_t\{\mathcal K(t)\}=1$, the norm of $T$ is $1$. Finally, given $f\in K$ choose $k$ such that $\|f-f_k\|_{APH_\Gamma(T)}\le\frac{\varepsilon}{3}$. Then we have by \eqref{eq3.1}
\begin{equation*}
\|Tf-f\|_{APH_\Gamma(T)} \leq \|T(f-f_k)\|_{APH_\Gamma(T)}+\|Tf_k-f_k\|_{APH_\Gamma(T)}+\|f_k-f\|_{APH_\Gamma(T)}<\varepsilon .
\end{equation*}
This completes the proof of the theorem.
\end{proof}

\subsection{Proof of Theorem \ref{totalspecprop}}

(1) The result follows directly from Remark \ref{equivway}.


\noindent (2) The fact that for $f\in SAP_\Sigma(S)$ the functions $\tilde{f}_k$ in Definition \ref{def1} can be chosen from $AP_{\Sigma(s,k)}(\mathbb R)$ follows from Theorem \ref{propconeq} and the definition of spectra of elements of $SAP(\partial \mathbb D)$. 

Let us show the validity of the converse statement.
Let $f\in SAP(S)$ and $s\in S$. Assume that for any $\varepsilon>0$ 
the functions $\tilde{f}_k$ in Definition \ref{def1} (for $f$ and $s$) can be chosen in $AP_{\Sigma(s,k)}(\mathbb R)$. 
Let $\rho$ be a smooth cut-off function equals $1$ in a neighbourhood of $s$ containing in the open set $\gamma_{s}^{-1}(\delta)\cup\{s\}\cup\gamma_{s}^1(\delta)$ and $0$ outside of this set. Let us consider a function 
$\hat f$ on $\partial\mathbb D\setminus\{s\}$ that coincides with $\rho f_1$ on $\gamma_{s}^1(\delta)$ and with
$\rho f_{-1}$ on $\gamma_{s}^{-1}(\delta)$ and equals $0$ outside of these arcs.
By the definition of spectra of functions in $SAP(\partial \mathbb D)$ the function $\hat f$ belongs to $SAP_{\Sigma}(S)$. Next, Theorem \ref{totalspecprop} (1) implies that
\begin{equation*}
|a_\lambda^k(\hat f,s)-a_\lambda^k(f,s)|=|a_\lambda^k(\hat f,s)-a_\lambda^k(\rho f,s)|<\varepsilon,\ \ \ k\in\{-1,1\}.
\end{equation*}
Since $\varepsilon>0$ is arbitrary, the latter inequality shows that if $a_\lambda^k(f,s) \ne 0$, then $\lambda \in \Sigma(s,k)$, as required.

\noindent (3) First, let us show that the set $T(f)$ of points of discontinuity of a function $f \in SAP(\partial \mathbb D)$ is at most countable. For each $s \in \partial \mathbb D$ define
\begin{equation*}
c_{s}(f):=\underset{\varepsilon \to 0+}{\overline{\lim}} \left( \sup_{s',s'' \in \bigl(e^{i(t-\varepsilon)},e^{i(t+\varepsilon)}\bigr)} |f(s')-f(s'')|  \right),
\end{equation*}
where $s=e^{it}$.
One has $c_s(f) \ne 0$ if and only if $s \in T(f)$. For $n\in\mathbb N$, we put $T_n(f):=\{s \in \partial \mathbb D: c_s(f) \geq \frac 1n\}$, so that $T(f)=\bigcup_{n=1}^\infty T_n(f)$. Suppose that $T(f)$ is uncountable, then there exists $n\in\mathbb N$ such that $T_n(f)$ is infinite. Since $\partial \mathbb D$ is compact, $T_{n}(f)$ has a limit point $e^{it_0}$. Choosing $\varepsilon<\frac{1}{2n}$ in Definition \ref{def1} (for $f$ and $e^{it_0}$) from the fact that $e^{it_0}$ is a limit point of $T_n(f)$ one obtains that the required functions $f_k$ do not exist, a contradiction. 

According to statement (2) $f_{k,s} \equiv \const$ for all points $s \in \partial \mathbb D$ at which $f$ is continuous. Therefore due to the previous statement $\speca_s^k(f)$ is $\{0\}$ or $\emptyset$ for all but at most countably many values of $s\in S$. Since for each $s \in S$ the spectrum $\speca_s^k(f)$ is at most countable, the required result follows.\hspace*{110mm} 

\subsection{Proof of Proposition \ref{impthm}}
\label{impthmsect}

Let $g \in APH(T)$, put $g_1(x):=g(x)$, $g_2(x):=g(x+i\pi)$, $x \in \mathbb R$. 
It follows easily from the approximation result for algebra $APH(T)$ cited in Section 3.1 that
$\speca(g_1)=\speca(g_2)$ and for each $\lambda \in \speca(g_1)$
\begin{equation*}
a_{\lambda}(g_1)=e^{\lambda \pi}a_{\lambda}(g_2).
\end{equation*}
Suppose that $f \in SAP(S) \cap H^\infty(\mathbb D)$.  Then Theorem \ref{maxidthm} (1) implies that for each $s \in S$ and $k\in\{-1,1\}$ the functions $f_{k,s}$ are the boundary values of the function $(i_\Sigma^s\circ\iota_{\Sigma(s)})^*(f)\in APH(T)$. 

This gives the required result.

\subsection{Proofs of Proposition \ref{strprop} and Theorem \ref{Bapprox}}
\label{Bapproxsect}

Our proof of Theorem \ref{Bapprox} is based on the equivalence established in Theorem \ref{equivthm}. We first formulate the $B$-valued analogues of the definitions of almost periodic and semi-almost periodic functions for $B$ being a complex Banach space.

Let $C^B_b(\mathbb R)$ and $C^B_b(T)$ denote the Banach spaces of $B$-valued bounded continuous functions on $\mathbb R$ and $T$, respectively, with norms $\|f\|:=\sup_{x}\|f(x)\|_B$.

\begin{definition}
\label{def3}
1) A function $f\in C^B_b(\mathbb R)$ is said to be \textit{almost periodic} if the family of its translates $\{S_\tau f\}_{\tau \in \mathbb R}$\,, $S_\tau f(x):=f(x+\tau)$, $x\in\mathbb R$, is relatively compact in $C^B_b(\mathbb R)$. 

2) 
A function $f\in C^B_b(T)$ is called \textit{holomorphic almost periodic} if it is holomorphic in the interior of $T$
and the family of its translates $\{S_x f\}_{x \in \mathbb R}$ is relatively compact in $C^B_b(T)$. 
\end{definition}

Let $AP^B(\mathbb R)$ and $APH^B(T)$ denote the Banach spaces of almost periodic and holomorphic almost periodic functions on $\mathbb R$ and $T$, respectively, endowed with $\sup$-norms. 

\begin{remark}
Since $AP(\mathbb R)$ and $APH(T)$ have the approximation property, it follows from Theorem \ref{equivthm}
that the functions of the form $\sum_{l=1}^m b_l e^{i\lambda_l x}$ ($x \in \mathbb R, b_l\in B, \lambda_l\in\mathbb R$) and 
$\sum_{l=1}^m b_l e^{i\lambda_l z}$ ($z \in T, b_l\in B, \lambda_l\in\mathbb R$) are dense in $AP^B(\mathbb R)$ and $APH^B(T)$, respectively. 
\end{remark}

As in the case of scalar almost periodic functions, a Banach-valued 
almost periodic function $f  \in AP_B(\mathbb R)$ is characterized
by its Bohr-Fourier coefficients $a_\lambda(f)$ and 
the spectrum $\speca(f)$, defined in terms of the mean value
\begin{equation}                   
\label{Bmeanf}
M(f):=\lim_{T \to +\infty}\frac{1}{2T}\int_{-T}^T f(x)dx \in B.
\end{equation}
Namely, we define
\begin{equation*}
a_\lambda(f):=M(fe^{-i \lambda x}), \quad \lambda \in \mathbb R.
\end{equation*}
It follows from the above remark and the properties of scalar almost periodic functions that 
$a_\lambda(f) \ne 0$ in $B$ for at 
most countably many values of $\lambda$.
These values constitute the spectrum $\speca(f)$ of $f$. E.g.,
if $f=\sum_{l=1}^\infty b_le^{i\lambda_l x}$ ($\sum_{l=1}^\infty\|b_l\|<\infty$ and all $b_l \ne 0$), then $\speca(f)=\{\lambda_1,\lambda_2,\dots\}$.

Similarly one defines the mean-values and the spectra for functions from $APH^B(T)$.


Let $L_B^\infty(\partial\mathbb D)$ be the Banach space of $B$-valued bounded measurable functions on $\partial\mathbb D$ equipped with $\sup$-norm.

\begin{definition}
\label{def22}
A function $f \in L_B^\infty(\partial\mathbb D)$ is called {\em semi-almost periodic} if for any $s \in \partial \mathbb D$ and any $\varepsilon>0$ there exist a number $\delta=\delta(s,\varepsilon) \in (0,\pi)$ and functions $f_k:\gamma_{s}^k(\delta)\to B$, $\gamma_{s}^k(\delta):=\{se^{ikx} :\ 0 \leq x<\delta<2\pi\}$, $k\in\{-1,1\}$, such that functions
\begin{equation*}
x \mapsto f_k\bigl(se^{ik\delta e^x}\bigr), \quad -\infty<x<0, \quad k \in \{-1,1\},
\end{equation*}
are restrictions of $B$-valued almost periodic functions from $AP^B(\mathbb R)$
and
\begin{equation*}
\sup_{z \in \gamma_{s}^k(\delta)}\|f(z)-f_k(z)\|_B<\varepsilon,\ \ \ k\in\{-1,1\}.
\end{equation*}
\end{definition}

Analogously to the scalar case, for a closed subset $S\subset\partial\mathbb D$ by $SAP^B(S) \subset L_B^\infty(\partial\mathbb D)$ we denote the Banach space of semi-almost periodic functions that are continuous on $\partial\mathbb D \setminus S$, so that $SAP(S):=SAP^{\mathbb C}(S)$.

Let $SAP(S) \otimes B$ denote the completion in $L_B^\infty(\partial\mathbb D)$ of the symmetric tensor product of $SAP(S)$ and $B$.

\begin{proposition}
\label{reprprop3}
$SAP^B(S)=SAP(S) \otimes B$.
\end{proposition}

The statement is an immediate consequence of Theorem \ref{equivthm} and the following two facts: each function $f\in SAP^B(S)$ admits a norm preserving extension to the maximal ideal space $b^S(\partial\mathbb D)$ of the algebra $SAP(S)$ as a continuous $B$-valued function, and $C(b^S(\partial\mathbb D))$ has the approximation property. The first fact follows straightforwardly from the definitions of $SAP^B(S)$ and $b^S(\partial\mathbb D)$ (see \cite{BK}) and the existence of analogous extensions of functions in $AP^B(\mathbb R)$ to $b\mathbb R$, while the second fact is valid for any algebra $C(X)$ on a compact Hausdorff topological space $X$ (it can be proved using finite partitions of unity of $X$). 

Next, we introduce the Banach space $AP^B(\mathbb R\sqcup (\mathbb R +i\pi)):=AP(\mathbb R\sqcup (\mathbb R +i\pi))\otimes B$ of $B$-valued almost periodic functions on $\mathbb R\sqcup (\mathbb R +i\pi)$, see Example \ref{basicex}. Also, for each $s\in S$ we define a linear isometry
$L_{s}^B:AP^B(\mathbb R\sqcup (\mathbb R +i\pi))\to L_B^\infty(\partial\mathbb D)$ by the formula (cf. \eqref{formF})
\begin{equation}
\label{formF1}
(L_{s}^B g)(z):=
(g\circ\Log\circ \varphi_{s})(z),\quad g\in AP^B(\mathbb R\sqcup (\mathbb R +i\pi)).
\end{equation}

Now, using Proposition \ref{reprprop3} we prove a $B$-valued analog of Theorem \ref{propconeq}.
\begin{theorem}
\label{Bpropconeq}
For every point $s\in\partial\mathbb D$ there exists a bounded linear operator $E_{s}^B: 
SAP^B(\partial\mathbb D)\to AP^B(\mathbb R\sqcup (\mathbb R +i\pi))$ of norm $1$ such that for each $f\in SAP^B(\partial\mathbb D)$ the function $f-L_{s}^B(E_{s}^B f)\in SAP^B(\partial\mathbb D)$ is continuous and equal to $0$ at $s$. Moreover, any bounded linear operator $SAP^B(\partial\mathbb D)\to AP^B(\mathbb R\sqcup (\mathbb R +i\pi))$ satisfying this property coincides with $E_{s}^B$.
\end{theorem}
\begin{proof}
According to Proposition \ref{reprprop3}, it suffices to define the required operator $E_{s}^B$ 
on the space of functions of the form $f=\sum_{l=1}^m b_l f_l$, where $b_l \in B$, $f_l \in SAP(\partial\mathbb D)$. In this case we set
\begin{equation*}
E_{s}^B(f):=\sum_{l=1}^m b_l E_{s}(f_l),
\end{equation*}
where $E_{s}$ is the operator from Theorem \ref{propconeq}.
Let ${\mathcal B_1}$ denote the unit ball in $B^*$. Then according to Theorem \ref{propconeq} we have
\begin{multline}
\notag
\|E_{s}^B(f)\|=\sup_{z \in \partial\mathbb D}\|E_{s}^B(f)(z)\|_B=
\sup_{z \in \partial\mathbb D,~\varphi \in {\mathcal B_1}}|\varphi(E_{s}^B(f)(z))|= \\
\sup_{z \in \partial\mathbb D,~\varphi \in {\mathcal B_1}} \left|\sum_{l=1}^m \varphi(b_l) E_{s}(f_l)(z) \right| \leq \sup_{z\in\partial\mathbb D,~\varphi \in {\mathcal B_1}}\left| \sum_{l=1}^m \varphi(b_l) f_l(z)\right|=\|f\|.
\end{multline}
This implies that $E_{s}^B$ is continuous and of norm $1$ on a dense subspace of $SAP^B(\partial\mathbb D)$. Moreover, for any function $f$ from this subspace we have (by Theorem \ref{propconeq}),
$f-L_{s}^B(E_{s}^B(f))\in SAP^B(\partial\mathbb D)$ is continuous and equal to $0$ at $s$. Extending $E_{s}^B$ by continuity to $SAP^B(\partial\mathbb D)$ we obtain the operator satisfying the required properties. Its uniqueness follows from the uniqueness of operator $E_{s}$.
\end{proof}

We make use of the functions $f_{-1,s}^B(x):=(E_{s}^B f)(x)$ and $f_{1,s}^B(x):=(E_{s}^B f)(x+i\pi)$, $x\in\mathbb R$, belonging to $AP_B(\mathbb R)$ to define the
left ($k=-1$) and the right ($k=1$) mean values of $f\in SAP^B(\partial\mathbb D)$ over $s$:
$$M_{s}^k(f):=M(f_{k,s}) \in B.$$
Then using formulas similar to those of the scalar case we define the Bohr-Fourier coefficients $a_\lambda^k(f,s) \in B$ and the spectrum $\speca_{s}^k(f)$ of $f$ over $s$.
It follows straightforwardly from the properties of the spectrum of 
a $B$-valued almost periodic function on $\mathbb R$ that $\speca_{s}^k(f)$ is at most countable. 

By $SAP^B_{\,\Sigma}(S) \subset SAP^B(S)$ we denote the Banach algebra of semi-almost periodic functions $f$ with
$\speca_s^k(f) \subset \Sigma(s,k)$
for all $s \in S$, $k \in \{-1,1\}$.
Note that $SAP_\Sigma^B(S)=SAP_\Sigma(S) \otimes B$, i.e.,
this definition is equivalent to the one used in Section 3.2 
(the proof is obtained easily from Definition \ref{def22}, using an appropriate partition of unity on $\partial \mathbb D$ and Theorems \ref{equivthm} and \ref{prop3.2}, see \cite{BK} for similar arguments).

Also, a statement analogous to Theorem \ref{impthm} holds for  $SAP_{\Sigma}^B(S)\cap H_B^\infty(\mathbb D)$. Namely, if $f \in SAP_{\Sigma}^B(S)\cap H_B^\infty(\mathbb D)$, then $$\speca_s^1(f)=\speca_s^{-1}(f)=:\speca_s(f).$$

\begin{proof}[Proof of Proposition \ref{strprop}]
We must show that if $f\in (SAP_{\Sigma}(S)\cap H^\infty(\mathbb D))_B$ on the maximal ideal space $b^S(\mathbb D)$ of algebra $SAP(S)\cap H^{\infty}(\mathbb D)$, then $f|_{\mathbb D}\in SAP_{\Sigma}^B(S)\cap H_B^\infty(\mathbb D)$.

Indeed, since $f\in C^B(b^S(\mathbb D))$ and $C(b^S(\mathbb D))$ has the approximation property, $f\in C(b^S(\mathbb D))\otimes B$ by Theorem \ref{equivthm}.
Next, $C(b^S(\mathbb D))$ is generated by algebra $SAP(S)\cap H^\infty(\mathbb D)$ and its conjugate. Therefore $f$ can be uniformly approximated on $b^S(\mathbb D)$ by a sequence of $B$-valued polynomials in variables from algebras $SAP(S)\cap H^\infty(\mathbb D)$ and its conjugate.
This easily implies that $f|_{\partial\mathbb D}$ is well defined and belongs to $SAP^B(S)$. In fact, $f|_{\partial\mathbb D}\in SAP_{\Sigma}^B(S)$ because $\phi(f)\in SAP_{\Sigma}(S)$ and the Bohr-Fourier coefficients of $f$ satisfy $a_{\lambda}^k(\varphi(f),s)=\varphi\bigl( a_{\lambda}^k(f,s)\bigr)$ for any $s\in S$, $k\in\{-1,1\}$ and $\phi\in B^*$. Further, by the definition
$f|_{\mathbb D}$ is such that 
$\varphi(f)\in H^\infty(\mathbb D)$ for any $\varphi \in B^*$.
This shows that $f\in H_B^\infty(\mathbb D)$.
\end{proof}

For the proof of Theorem \ref{Bapprox} we require some auxiliary results. 
Let $APC(T)$ be the Banach algebra of functions $f:T \to \mathbb C$ uniformly continuous on $T$ and almost periodic on each horizontal line. 
We define $APC^B(T):=APC(T) \otimes B$. The proof of the next statement is analogous to the proof of Lemma 4.3 in \cite{BK}.

\begin{lemma}
\label{fivelem}
Suppose that $f_1 \in AP^B(\mathbb R)$, $f_2 \in AP^B(\mathbb R+i\pi)$. Then there exists a function $F \in APC^B(T)$ harmonic in the interior of $\Sigma$ whose boundary values are $f_1$ and $f_2$. Moreover, $F$ admits a continuous extension to the maximal ideal space $bT$ of $APH(T)$.
\end{lemma}

The proof of the next statement uses Lemma \ref{fivelem} and is very similar to the proof of Lemma 4.2 (for $B=\mathbb C$) in \cite{BK}, so we omit it as well.

Suppose that $s:=e^{it}$ and $\gamma_{s}^{k}(\delta)$ are arcs defined in \eqref{arcs}. For $\delta\in (0,\pi)$ we set $\gamma_{1}(s,\delta):={\rm Log}(\varphi_{s}(\gamma_{s}^{1}(\delta)))\subset\mathbb R$ and 
$\gamma_{-1}(s,\delta):={\rm Log}(\varphi_{s}(\gamma_{s}^{-1}(\delta)))\subset \mathbb R+i\pi$ (see Example \ref{ex2}). 

\begin{lemma}
\label{fourlem}
Let $s \in S$. Suppose that $f \in SAP^B(\{-s,s\})$. We put $f_k=f|_{\gamma_s^k(\pi)}$ and 
define on $\gamma_{k}(s,\pi)$ 
$$h_k:=f_k \circ \varphi_{s}^{-1} \circ \Log^{-1}, \quad k \in \{-1,1\}.$$ 
Then for any $\varepsilon>0$ there exist $\delta_\varepsilon \in (0,\pi)$ and a 
function $H\in APC^B(T)$ harmonic in the interior $T_0$ of $T$ such that 
\begin{equation*}
\sup_{z \in \gamma_{k}(s,\delta_\varepsilon)}\|h_k(z)-H(z)\|_B<\varepsilon, \quad k \in \{-1,1\}.
\end{equation*}
\end{lemma}

Let $s \in \partial\mathbb D$ and $U_{s}$ be the intersection of an open disk of some radius $ \leq 1$
centered at $s$ with $\bar{\mathbb D} \setminus s$. We call such $U_{s}$ a {\em circular neighbourhood} of $s$.

\begin{definition}
\label{circdef2}
We say that a bounded continuous function $f:\mathbb D \to B$ is {\em almost-periodic near} $s$ if there exist a circular neighbourhood $U_{s}$, and a function $\bar{f} \in APC^B(T)$ such that
\begin{equation}
\label{circdef}
f(z)=(L_{s}^B\bar{f})(z):=\bigl(\bar{f} \circ \Log \circ \varphi_{s}\bigr)(z), \quad z \in U_{s}.
\end{equation}
\end{definition}

In what follows for $\Sigma:S\times\{-1,1\}\to 2^{\mathbb R}$ such that $\Sigma(s)=\mathbb R$ for each $s\in S$ we omit writing $\Sigma$ in  $a_\Sigma^S$, $i_\Sigma^s$, $b_\Sigma^S(\mathbb D)$ etc., see Section~3.4 for the corresponding definitions.

In the proof of Theorem 1.8 of \cite{BK} (see Lemmas 4.4, 4.6 there) we established, cf. Theorem \ref{maxidthm},
\begin{itemize}
\item[(1)]
Any scalar harmonic function $f$ on $\mathbb D$ almost periodic near $s$ admits a continuous extension $f_{s}$  to $(a^S)^{-1}(\bar{U}_{s})\subset b^{S}(\mathbb D)$ for some circular neighbourhood $U_{s}$.

\item[(2)]
For any $s \in S$ and any $g \in APH(T)$  the holomorphic function $\tilde g:=L_{s}g$ on $\mathbb D$ almost periodic near $s$ is such that $\tilde g_{s}\circ i^{s}$ coincides with the extension of $g$ to $bT$. 
\end{itemize}


More generally, Lemma \ref{fivelem}, statements (1) and (2) and the fact that $AP^B(\mathbb R)=AP(\mathbb R) \otimes B$ (see Theorems \ref{equivthm} and \ref{prop3.1}) imply

\begin{itemize}
\item[(3)]
Any $B$-valued harmonic function $f$ on $\mathbb D$ almost periodic near $s$ admits a continuous extension $f_{s}^B$ to $(a^S)^{-1}(\bar{U}_{s})\subset b^{S}(\mathbb D)$ for some circular neighbourhood $U_{s}$.

\item[(4)]
For any $s \in S$ and any $g \in APH^B(T)$
the $B$-valued holomorphic function $\tilde g:=L_{s}^B g$ on $\mathbb D$ almost periodic near $s$ is such that $\tilde g_{s}^B\circ i^{s}$ coincides with the extension of $g$ to $bT$. 
\end{itemize}

\begin{lemma}\label{l5.6}
Let $f\in SAP^B_{\Sigma}(S)\cap H_B^{\infty}(\mathbb D)$ and $s\in\partial\mathbb D$. There is a bounded $B$-valued holomorphic function $\hat f$ on $\mathbb D$ almost periodic near $s$ such that for any $\varepsilon>0$
there is a circular neighbourhood $U_{s;\varepsilon}$ of $s$ so that
$$
\sup_{z\in U_{s;\varepsilon}}||f(z)-\hat f(z)||_B<\varepsilon.
$$
Moreover, $\hat{f}=L_{s}^B\bar{f}$ for some $\bar{f} \in APH^B_{\Sigma(s)}(T)$.
\end{lemma}

\begin{proof}
Assume, first, that $s\in S$. By Lemma \ref{fourlem}, for any $n\in\mathbb N$ there exist a number $\delta_n \in (0,\pi)$ and a function $H_n\in APC^B(T)$ harmonic on $T_0$ such that 
\begin{equation}\label{e5.3}
\sup_{z \in \gamma_k(s,\delta_n)}\|f_k(z)-H_n(z)\|_B<\frac{1}{n}\, , \quad k \in \{-1,1\}.
\end{equation}
Using the Poisson integral formula for the bounded $B$-valued harmonic function $f-h_n$,
$h_n:=L_{s}^B H_n:= H_n\circ \Log\circ\varphi_{s}$,
on $\mathbb D$ we easily obtain from \eqref{e5.3} that there is a circular neighbourhood $V_{s;n}$ of $s$ such that
\begin{equation}\label{e5.4}
\sup_{z \in V_{s;n}}\|f(z)-h_n(z)\|_B<\frac{2}{n} .
\end{equation}
According to (3) each $h_n$ admits a continuous extension $\hat h_n$ to $(a^S)^{-1}(s)\cong bT$. Moreover, \eqref{e5.4} implies that the restriction of the sequence $\{\hat h_n\}_{n\in\mathbb N}$ to $(a^S)^{-1}(s)$ forms a Cauchy sequence in $C^B((a^S)^{-1}(s))$. Let $\hat h\in C^B((a^S)^{-1}(s))$ be the limit of this sequence. 

Further, for any functional $\phi\in B^*$ the 
function $\phi\circ f\in SAP(S)\cap H^\infty(\mathbb D)$ 
and therefore admits a continuous extension $f_{\phi}$ 
to $(a^S)^{-1}(s)$ such that on $(i^{s})^{-1}\left((a^S)^{-1}(s)\right)$ the function $f_{\phi}\circ i^{s}$ belongs to 
$APH(T)$. Now, \eqref{e5.4} implies 
that $f_\phi=\phi\circ\hat h$ for any $\phi\in B^*$. Then it follows from Theorems \ref{equivthm} and \ref{prop3.2}
that $\hat h\circ i^{s}\in APH^B(T)$. Therefore by (4) we 
find a bounded $B$-valued holomorphic function $\hat f$ on $\mathbb D$ 
of the same $\sup$-norm as $\hat h$ almost periodic near $s$ such 
that its extension to $(a^S)^{-1}(s)$ coincides with $\hat h$. 
Next, by the definition of the topology of $b^S(\mathbb D)$, see \cite{BK}, Lemma 4.4 (a), we obtain that for any $\varepsilon>0$ there is a number $N\in\mathbb N$ such that for all $n\geq N$,
$$
\sup_{z \in V_{s;n}}\|\hat f(z)-h_n(z)\|_B<\frac{\varepsilon}{2}.
$$ 
Now, choose $n\geq N$ in \eqref{e5.4} such that the right-hand side there is $<\frac{\varepsilon}{2}$. For this $n$ we set $U_{s;\varepsilon}:=V_{s;n}$. Then the previous inequality and \eqref{e5.4} imply the required inequality
$$
\sup_{z\in U_{s;\varepsilon}}||f(z)-\hat f(z)||_B<\varepsilon .
$$

Further, if $s\not\in S$, then, by definition, $f|_{\partial\mathbb D}$ is continuous at $s$. In this case as the function $\hat f$ we can choose the constant $B$-valued function equal to $f(s)$ on $\mathbb D$. Then the required result follows from the Poisson integral formula for $f-\hat f$.

By definition, $\hat{f}$ is determined by formula (\ref{circdef}) with an $\bar{f} \in APH^B(T)$. Let us show that $\bar{f} \in APH_{\Sigma(s)}^B(T)$. To this end it suffices to prove that $\varphi(\bar{f}) \in APH_{\Sigma(s)}(T)$ for any $\varphi \in B^*$. Indeed, it follows from the last inequality that the extension of $\varphi(f) \in SAP_\Sigma(S) \cap H^\infty(\mathbb D)$ to $(a^S_\Sigma)^{-1}(s)$ coincides with $\varphi(\bar{f})$. By the definition of spectrum of a semi-almost periodic function, this implies that $\speca_{s}(\varphi(\bar{f})) \subset \Sigma(s)$.
\end{proof}

Now, we are ready to prove Theorem \ref{Bapprox}.

\vspace*{2mm}

The inclusion $A_\Sigma^S \subset SAP_\Sigma(S) \cap H^\infty(\mathbb D)$ follows from Example \ref{ex5}. Indeed, for $s \in S$ assume that
the holomorphic function $e^{\lambda h}\in H^\infty(\mathbb D)$ is such that $\Real(h)|_{\partial\mathbb D}$ is the characteristic function of the closed arc going in the counterclockwise direction from the initial point at $s$ to the endpoint at $-s$ and such that $\frac{\lambda}{\pi}\in\Sigma(s)$. Then Example \ref{ex5} implies that $e^{\lambda h}\in SAP(\{s,-s\})\cap H^\infty(\mathbb D)$ and $\speca_{s}(e^{\lambda h})=\{\frac{\lambda}{\pi}\}$. In particular, 
$(z+s)e^{\frac{\lambda}{\pi}h} \in SAP_{\Sigma|_{\{s\}}}(\{s\}) \cap H^\infty(\mathbb D)$, as required.

\vspace*{2mm}

Let us prove the opposite inclusion.

\vspace*{2mm}

(A) Consider first the case $S=F$, where $F=\{s_i\}_{i=1}^m$ is a finite subset of $\partial\mathbb D$. 
Let $f \in SAP^B_\Sigma(F) \cap H^\infty_B(\mathbb D)$. Then according to Lemma \ref{l5.6} there exists a function $f_{s_1} \in APH^B_{\Sigma(s_1)}(T)$
such that the bounded $B$-valued holomorphic function $g_{s_1}-f$, where $g_{s_1}:=f_{s_1}\circ\Log\circ\varphi_{s_1}$, on $\mathbb D$ is continuous and equals $0$ at $s_1$.

Let us show that $g_{s_1} \in A^{\{s_1,-s_1\}} \otimes B$. 
Since $f_{s_1} \in APH^B_{\Sigma(s_1)}(T)$, by Theorems \ref{equivthm} and \ref{prop3.2} it can be approximated in $APH^B_{\Sigma(s_1)}(T)$ by finite sums of functions of the form $b e^{i\lambda z}$, $b \in B$, $\lambda \in \Sigma(s_1)$, $z \in T$. In turn, $g_{s_1}$ can be approximated by finite sums of functions of the form $b e^{i\lambda \Log \circ \varphi_{s_1}}$, $b\in B$. 
As was shown in \cite{BK}, $e^{i\lambda \Log \circ \varphi_{s_1}} \in A^{\{s_1,-s_1\}}$. Hence, $g_{s_1} \in A^{\{s_1,-s_1\}} \otimes B$.
We define
\begin{equation*}
\hat{g}_{s_1}=\frac{g_{s_1}(z)(z+s_1)}{2s_1}.
\end{equation*}
Then, since the function $z \mapsto (z+s_1)/(2s_1) \in A(\mathbb D)$ and equals $0$ at $-s_1$, and $g_{s_1} \in A^{\{s_1,-s_1\}} \otimes B$, the function $\hat{g}_{s_1} \in A^{\{s_1\}} \otimes B$. Moreover, by the construction of $\hat{g}_1$ and the definition of the spectrum $\hat{g}_{s_1} \in A_{\Sigma(s_1)}^{\{s_1\}} \otimes B$. By definition, the difference $\hat{g}_{s_1}-f$ is continuous and equal to zero at $z_1$.  Thus, 
$$
\hat{g}_{s_1}-f \in SAP^B_{\Sigma|_{F \setminus \{s_1\}}}(F \setminus \{s_1\})\cap H_B^\infty(\mathbb D).
$$
We proceed in this way to get functions $\hat{g}_{s_k} \in A_{\Sigma(s_k)}^{\{s_k\}} \otimes B$, $1\le k\le m$, such that 
\begin{equation*}
f -\sum_{k=1}^m \hat{g}_{s_k} \in A^B(\mathbb D),
\end{equation*}
where $A^B(\mathbb D)$ is the Banach space of $B$-valued bounded holomorphic functions on $\mathbb D$ continuous up to the boundary. As in the scalar case using the Taylor expansion at $0$ of functions from $A^B(\mathbb D)$ one can easily show that $A^B(\mathbb D)=A(\mathbb D)\otimes B$. Therefore, $f\in A^F_{\Sigma}\otimes B$.

\vspace*{2mm}

(B) Let us consider the general case of $S \subset \partial\mathbb D$ being an arbitrary closed set. 
Let $f \in SAP^B_\Sigma(S) \cap H^\infty_B(\mathbb D)$. As follows from Lemma \ref{l5.6} and the arguments presented in part (A), given an $\varepsilon>0$ there exist points $s_k \in \partial\mathbb D$, functions $f_k \in A_{\Sigma(s_k)}^{\{s_k\}} \otimes B$ and circular neighbourhoods $U_{s_k}$ ($1 \leq k \leq m$) such that $\{U_{s_k}\}_{k=1}^m$ forms an open cover of $\partial\mathbb D \setminus \{s_k\}_{k=1}^m$ and
\begin{equation}
\label{b3}
\|f(z)-f_k(z)\|_B<\varepsilon\ \ \text{ on }\ \ U_{s_k}, \quad 1 \leq k \leq m.
\end{equation}
Since $S$ is closed, for $s_k \not\in S$ we may assume that $f_k$ is continuous in $\bar{U}_{s_k}$.
Let us define a $B$-valued $1$-cocycle $\{c_{kj}\}_{k,j=1}^m$ on intersections of the sets in $\{U_{s_k}\}_{k=1}^m$ by the formula
\begin{equation}\label{e5.6'}
c_{kj}(z):=f_k(z)-f_j(z), \quad z \in U_{s_k} \cap U_{s_j}.
\end{equation}
Then \eqref{b3} implies
$\sup_{k,j,z}||c_{kj}(z)||_B<2\varepsilon$.
Let $A\Subset\cup_{k=1}^m U_{s_k}$ be an open annulus with outer boundary $\partial\mathbb D$. Using the argument from the proof of Lemma 4.7 in \cite{BK} one obtains that if the width of the annulus is sufficiently small, then there exist $B$-valued functions $c_i$ holomorphic on $U_{s_i}\cap A$ and continuous on $\bar{U}_{s_i}\cap\bar{A}$  satisfying
\begin{equation}\label{e5.12'}
\sup_{z\in U_{s_i}\cap A}||c_i(z)||_B\leq 3\varepsilon
\end{equation}
and such that
\begin{equation}\label{e5.12}
c_i(z)-c_j(z)=c_{ij}(z), \quad z \in U_{s_i} \cap U_{s_j}\cap A.
\end{equation}
For such $A$ let us define a function $f_\varepsilon$ on $\bar{A} \setminus \{s_i\}_{i=1}^m$ by formulas
\begin{equation*}
f_\varepsilon(z):=f_i(z)-c_i(z), \quad z \in U_{s_i}\cap\bar{A}.
\end{equation*}
According to \eqref{e5.6'} and \eqref{e5.12}, $f_{\varepsilon}$ is a bounded continuous $B$-valued function on $\bar{A}\setminus \{s_i\}_{i=1}^m$ holomorphic in $A$.
Furthermore, since $c_i$ is continuous on $\bar{U}_{s_i}\cap\bar{A}$, and $f_i\in A_{\Sigma(s_i)}^{\{s_i\}} \otimes B$ for $s_i\in S$, and $f_i\in A^B(\mathbb D)$ otherwise,  $f_\varepsilon|_{\partial\mathbb D}\in SAP^B_{\Sigma|_F}(F)$, where $F=\{s_i\}_{i=1}^m \cap S$. Also, from inequalities (\ref{b3}) and \eqref{e5.12'} we obtain
\begin{equation}
\label{ineq1t}
\sup_{z\in A}\|f(z)-f_\varepsilon(z)\|_B<4\varepsilon .
\end{equation}
Next, as in \cite{BK} we consider a $1$-cocylce subordinate to a cover of the unit disk $\mathbb D$ consisting of an open annulus having the same interior boundary as $A$ and the outer boundary $\{z\in\mathbb C\ :\ |z|=2\}$, and of an open disk centered at $0$ not containing $A$ but intersecting it by a nonempty set. 
Resolving this cocycle\footnote{There is a misprint in \cite{BK} at this place: instead of the inequality $\max_i\|\nabla\rho_i\|_{L^\infty(\mathbb C)}\le\widetilde C w(B\cap A)<\widetilde C\varepsilon$ for smooth radial functions $\rho_1,\rho_2$, it must be $\max_i\|\nabla\rho_i\|_{L^\infty(\mathbb C)}\le\frac{\widetilde C}{w(B\cap A)}$.} one obtains a $B$-valued holomorphic function $F_{\varepsilon}$ on $\mathbb D$ such that for an absolute constant $\hat{C}>0$
\begin{equation*}
\sup_{z\in\mathbb D}\|f(z)-F_\varepsilon(z)\|_B<\hat{C}\varepsilon 
\end{equation*}
and by definition $F_\varepsilon \in SAP_{\Sigma|_F}^B(F) \cap H^\infty_B(\mathbb D)$, where $F=\{s_1,\dots,s_m\} \cap S$.
 The latter inequality and part (A) of the proof show that the complex vector space generated by spaces $A_{\Sigma|_F}^F\otimes B$ for all possible finite subsets $F\subset S$ is dense in $SAP_{\Sigma}^B(S)\cap H_B^\infty(\mathbb D)$. Since by definition the closure of all such $A_{\Sigma|_F}^F \otimes B$ is
$A^S_{\Sigma} \otimes B$, we
obtain the required: $SAP_{\Sigma}^B(S)\cap H_B^\infty(\mathbb D)=A_\Sigma^S \otimes B$.

\subsection{Proof of Theorems \ref{isothm} and \ref{convthm}}

\begin{proof}[Proof of Theorem \ref{isothm}]
Corollary 1.6 in \cite{BK} states that $\kappa^*|_{\partial \mathbb D}:C(\partial \mathbb D) \mapsto C(\partial \mathbb D)$, the pullback by $\kappa|_{\partial \mathbb D}$, maps $SAP(\kappa_* S)$ isomorphically onto $SAP(S)$.
Following closely the arguments in its proof, one obtains even more: $\kappa^*$ maps $SAP_{\kappa_*\Sigma}(\kappa_* S)$ isomorphically onto $SAP_{\Sigma}(S)$. Since $\kappa^*$ preserves $H^\infty(\mathbb D)$, the required result follows.
\end{proof}

\begin{proof}[Proof of Theorem \ref{convthm}]
Let $f \in SAP(\partial \mathbb D) \cap H^\infty(\mathbb D)$. According to Lemma \ref{l5.6} there exists a function $f_{s}\in APH(T)$ such that the difference
\begin{equation*}
h:=f-F_{s},
\end{equation*}
where $F_{s}:=f_{s}\circ\Log\circ\varphi_{s}$, see \eqref{circdef}, is continuous and equal to $0$ at $s$. Therefore, it suffices to prove the assertion of the theorem for $F_{s}$. Let $\{z_n'\} \subset T_0$ and $\{s_n'\} \subset \mathbb R\cup (\mathbb R+i\pi)$ be the images of sequences $\{z_n\}$ and $\{s_n\}$ under the mapping $\Log \circ \varphi_{s}$ (see Example \ref{ex2}). By the hypotheses of the theorem we have $z_n'$, $s_n' \to \infty$ and $|z_n'-s_n'| \to 0$ as $n \to \infty$ (this follows from condition (\ref{convcond})). Since any function in $APH(T)$ is uniformly continuous (see Section 3), the latter implies the required result.
\end{proof}

\subsection{Proof of Theorems  \ref{coronathm} and \ref{holmapsthm}}

\begin{proof}[Proof of Theorem \ref{coronathm}]
In what follows we identify $\iota(\mathbb D)\subset b_\Sigma^S(\mathbb D)$ with $\mathbb D$, see Section~3.4.
By Theorem \ref{maxidthm}, the maximal ideal space $b_{\Sigma}^S(\mathbb D)$ is  $(\bar{\mathbb D} \setminus S)\sqcup\left(\sqcup_{s\in S}\,i_\Sigma^s(b_{\Sigma(s)}(T))\right)$ (here $i_\Sigma^s:b_{\Sigma(s)}(T)\to (a_\Sigma^S)^{-1}(s)$ is a homeomorphism).
For each $s \in S$ one has the natural map $\iota_{\Sigma(s)}:T \hookrightarrow b_{\Sigma(s)}(T)$ (determined by evaluations at points of $T$). Also,
the argument of the proof of Theorem 1.12 in \cite{BK} implies that the closure of $\mathbb D$ in $b_{\Sigma}^S(\mathbb D)$ contains (as a dense subset)
$(\bar{\mathbb D} \setminus S)\sqcup\left(\sqcup_{s\in S}\,i_\Sigma^s(\iota_{\Sigma(s)}(T))\right)$. Thus in order to prove the theorem, it suffices to show that $\iota_{\Sigma(s)}(T)$ is dense in $b_{\Sigma(s)}(T)$ if and only if $\Sigma(s)$ is a group. 

We will use the following result.

\begin{theorem}[\cite{RS}]
\label{rsthm}
Suppose that $\Gamma$ is the intersection of an additive subgroup of $\mathbb R$ and $\mathbb R_+$. Then the image of the upper half-plane $\mathbb H^+$ in the maximal ideal space $b_\Gamma(T)$ is dense.
\end{theorem}
Observe that in this case each element of $APH_\Gamma(T)$ is extended to a holomorphic almost periodic function on $\mathbb H^+$ by means of the Poisson integral. Therefore the evaluations at points of $\mathbb H^+$ of the extended algebra determine the map $\mathbb H^+\to b_\Gamma(T)$ of the theorem.


First, assume that $\Sigma(s)$ is a group. We have to show that $\iota_{\Sigma(s)}(T)$ is dense in $b_{\Sigma(s)}(T)$. 
Suppose that this is wrong. Then there exists $\xi \in b_{\Sigma(s)}(T)$ and a neighbourhood of $\xi$
\begin{equation*}
U(\lambda_1,\dots,\lambda_m,\xi,\varepsilon):=\{\eta \in b_{\Sigma(s)}(T): |\eta(e^{i\lambda_k z})-c_k|<\varepsilon ,\ 1\le k\le m\},
\end{equation*}
where $\lambda_1,\dots,\lambda_m \in \Sigma(s)$, $c_k:=\xi(e^{i\lambda_k z})$,
such that
$U(\lambda_1,\dots,\lambda_m,\xi,\varepsilon) \cap \cl\left(\iota_{\Sigma(s)}(T)\right)=\varnothing$, cf. the proof of Theorem~2.4 in \cite{BK}.
Therefore, 
\begin{equation}
\label{geqid}
\max_{1 \leq k \leq m}|e^{i\lambda_k z}-c_k| \geq \varepsilon>0 \quad \text{ for all } z \in T.
\end{equation}
Without loss of generality we may assume that $c_k \ne 0$ and $\lambda_k>0$, i.e., $e^{i\lambda_k z}-c_k \in APH_{\Sigma(s) \cap \mathbb R_+}(T)$. (For otherwise we replace $e^{i\lambda_k z}-c_k$ with $e^{-i\lambda_k z}-c_k^{-1}$. Here $e^{-i\lambda_k z}-c_k^{-1} \in APH_{\Sigma(s) \cap \mathbb R_+}(T)$ since $\Sigma(s)$ is a group. Also, (\ref{geqid}) will be satisfied, possibly with a different $\varepsilon>0$.) Note that $e^{i\lambda_k z}-c_k$ is not invertible in $APH_{\Sigma(s)}(T)$, since $\xi(e^{i\lambda_k z}-c_k)=0$. Therefore, since each function $e^{i\lambda_k z}-c_k$ is periodic (with period $\frac{2\pi}{\lambda_k}$), it has a zero in $T$. Since solutions of the equation $e^{i\lambda_k z}=c_k$ are of the form
\begin{equation*}
z_k=-\frac{i\ln |c_k|}{\lambda_k}+\frac{\Arg c_k+2\pi l}{\lambda_k}, \quad l \in \mathbb Z,
\end{equation*}
all zeros of $e^{i\lambda_k z}-c_k$ belong to $T$. Hence, in virtue of inequality (\ref{geqid}), one has
\begin{equation*}
\max_{1 \leq k \leq m}|e^{i\lambda_k z}-c_k| \geq \tilde{\varepsilon}>0 \quad \text{ for all } z \in \mathbb H^+.
\end{equation*}
This implies, by Theorem \ref{rsthm}, that there exist functions $g_1,\dots,g_m \in APH_{\Sigma(s) \cap \mathbb R_+}(\mathbb H_+)$ such that
\begin{equation*}
\sum_{k=1}^m g_k(z) (e^{i\lambda_k z}-c_k) = 1 \quad \text{for all }\quad z\in\mathbb H^+. 
\end{equation*}
In particular, the above identity holds on $T$. This gives a contradiction with the assumption $\xi(e^{i\lambda_k z}-c_k)=0$, $1\le k\le m$.

Now, assume that $\Sigma(s)$ is not a group, i.e., it contains a non-invertible element $\lambda_0$. Suppose that $\iota_{\Sigma(s)}(T)$ is dense in $b_{\Sigma(s)}T$. Then, since the modulus of $f_1:=e^{i\lambda_0 z}$ is bounded from below on $T$ by a positive number, there exists $g_1 \in APH_{\Sigma(s)}(T)$ such that $f_1g_1 \equiv 1$. Therefore, $g_1=e^{-i\lambda_0 z}\in APH_{\Sigma(s)}(T)$, i.e., $-\lambda_0 \in \Sigma(s)$, a contradiction.
\end{proof}

\begin{proof}[Proof of Theorem \ref{holmapsthm}]
For the proof we will need the following auxiliary result. 

Let $\Gamma\subset\mathbb R$ be a nontrivial additive semi-group. For a subset $X\subset T$ by $X_{\infty}$ we denote the set of limit points of $\iota_{\Gamma}(X)$ in $b_{\Gamma}(T)\setminus \iota_\Gamma(T)$.

\begin{lemma}\label{l4.12}
Let $G \in C([0,1),T)$ be such that
the closure of $G([0,1))$ in $T$ is non-compact. Then the set 
$G_\infty$ contains more than one element.
\end{lemma}

\begin{proof}
If there exists a horizontal line $\mathbb R+ic$, $0 \leq c \leq \pi$,
such that $\dist_T\bigl(G(t),\mathbb R+ic\bigr) \to 0$ as $t \to 1-$, then clearly $(\mathbb R+ic)_\infty=G_\infty$. Moreover, $(\mathbb R+ic)_\infty$ is infinite (e.g., it contains a subset homeomorphic to interval $[0,1]$).
In the case that such a line does not exist, one can find two closed substrips $T_1$, $T_2 \subset T$, $T_1 \cap T_2=\varnothing$,
such that the closures in $T$ of both $G([0,1)) \cap T_1$ and  $G([0,1)) \cap T_2$ are non-compact. Then $(G([0,1)) \cap T_1)_{\infty}$ and $(G([0,1)) \cap T_2)_{\infty}$ are nonempty,
while $(T_1)_{\infty} \cap (T_2)_{\infty}=\varnothing$. This implies the required statement.
\end{proof}

Now, we are ready to prove the theorem.
Suppose on the contrary that for a continuous map $F: [0,1]\to b_{\Sigma}^S(\mathbb D)$ the conclusion of the theorem is not valid. First, assume that there exists a point $c\in [0,1)$ such that $F(c)\in\bar{\mathbb D}\setminus S$ but $F([0,1])\not\subset\bar{\mathbb D}\setminus S$. Then, because $b_{\Sigma}^S(\mathbb D)\setminus (\bar{\mathbb D}\setminus S)$ is a compact set (here we naturally identify $\bar{\mathbb D}\setminus S$ we a subset of $b_\Sigma^S(\mathbb D)$), passing to a subinterval, if necessary, we may assume without loss of generality that $F[0,1) \subset \bar{\mathbb D} \setminus S$ and $F(1) \in (a_\Sigma^S)^{-1}(s)$ for a certain $s\in S$.
Define 
$$
G(t):=(\Log \circ \varphi_{s})\bigl(F(t)\bigr) \subset T, \quad t \in [0,1)
$$ 
(cf. Example \ref{basicex}).
Then $G$ satisfies conditions of Lemma \ref{l4.12} for $\Gamma=\Sigma(s)$.
Next, consider an $f \in SAP_{\Sigma}(S)$. According to Lemma \ref{l5.6} there exists a (unique) function $f_{s}\in APH_{\Sigma(s)}(T)$ such that the difference
$f-F_{s}$,
where $F_{s}:=f_{s}\circ\Log\circ\varphi_{s}$, is continuous and equal to $0$ at $s$. This yields 
$$
\lim_{t\to 1-}\bigl(f_s(G(t))-f(F(t))\bigr)=0.
$$
The latter implies that the set of limit points of $F([0,1))$ in $b_{\Sigma}^S(\mathbb D)\setminus (\bar{\mathbb D}\setminus S)$ 
is in one-to-one correspondence with the set of limit points $G_\infty$ of $\iota_{\Sigma(s)}\left(G([0,1))\right)$ in $b_{\Sigma(s)}(T)\setminus \iota_{\Sigma(s)}(T)$. 
By our assumption the set of limit points of $F([0,1))$ in $b_{\Sigma}^S(\mathbb D)\setminus (\bar{\mathbb D}\setminus S)$ consists of the point $F(1)$. This contradicts the assertion of Lemma \ref{l4.12}. Hence, in this case $F([0,1])\subset\bar{\mathbb D}\setminus S$. 

In the second case, $F([0,1])\subset b_{\Sigma}^S(\mathbb D)\setminus (\bar{\mathbb D}\setminus S)$. Let $s\subset S$ be such that $F([0,1])\cap (a_{\Sigma}^S)^{-1}(s)\neq\emptyset$. Consider the continuous map
$\omega_{s}: b_{\Sigma}^S(\mathbb D)\to b_{\Sigma|_{\{s\}}}^{\{s\}}(\mathbb D)$ transpose to the embedding $SAP_{\Sigma|_{\{s\}}}(\{s\})\cap H^\infty(\mathbb D)\subset SAP_{\Sigma}(S)\cap H^\infty(\mathbb D)$. According to the case considered above, if $\omega_s\circ F:[0,1]\to b_{\Sigma|_{\{s\}}}^{\{s\}}(\mathbb D)$ is such that $(\omega_{s}\circ F)(c)\in\bar{\mathbb D}\setminus\{s\}$ for some $c\in [0,1)$, then $(\omega_s\circ F)([0,1])\subset \bar{\mathbb D}\setminus\{s\}$ which contradicts the assumption $F([0,1])\cap (a_{\Sigma}^S)^{-1}(s)\neq\emptyset$. Thus $(\omega_s\circ F)([0,1])\subset b_{\Sigma|_{\{s\}}}^{\{s\}}(\mathbb D)\setminus (\bar{\mathbb D}\setminus\{s\})=(a_{\Sigma|_{\{s\}}}^{\{s\}})^{-1}(s)$. This implies that $F([0,1])\subset \left(a_{\Sigma}^S\right)^{-1}(s)$.
\end{proof}

\subsection{Proof of Theorem \ref{shilovthm}}

Since $SAP_\Sigma(S) \cap H^\infty(\mathbb D)$ is generated by algebras $SAP_{\Sigma|_F}(F) \cap H^\infty(\mathbb D)$ for all possible finite subsets $F$ of $S$, the inverse limit of $\{K_{\Sigma|_{F}}^F\, ;\, \omega\}_{F\subset S\, ;\,\#F<\infty}$ of the corresponding \v{S}ilov boundaries  coincides with $K_\Sigma^S$ (see Section~3.4 for the corresponding notation). Therefore to establish the result it suffices to prove that
\begin{equation}
\label{shilov1}
K_{\Sigma|_{F}}^F=\left(\bigcup_{s\in F} i_{\Sigma|_{F}}^s\bigl(\cl_{\Sigma(s)}(\mathbb R) \cup \cl_{\Sigma(s)}(\mathbb R+i\pi)\bigr)\right) \cup \partial \mathbb D \setminus F.
\end{equation}

Since each point of $\partial\mathbb D\setminus F$ is a peak point for $A(\mathbb D)\, (\subset SAP_\Sigma|_{F}(F)\cap H^\infty(\mathbb D))$,  $\partial\mathbb D\setminus F\subset K_{\Sigma|_{F}}^F$.  Next, the closure of $\partial\mathbb D\setminus F$ in $b_{\Sigma|_{F}}^F(\mathbb D)$ (the maximal ideal space of $SAP_\Sigma|_{F}(F)\cap H^\infty(\mathbb D)$) coincides with the right-hand side of \eqref{shilov1}, see the proof of Theorem~1.14 in \cite{BK}. Thus the right-hand side of \eqref{shilov1} is a subset of 
$K_{\Sigma|_{F}}^F$. Finally, Theorem~1.14 of \cite{BK}  implies that for each $f\in SAP_\Sigma|_{F}(F)\cap H^\infty(\mathbb D)$,\, $|f|$ attains its maximum on the set in the right-hand side of \eqref{shilov1}. This produces the required identity.

One can easily show that the inverse limit of the family of sets in the right-hand sides of equations \eqref{shilov1} coincides
with $\left(\bigcup_{s\in S} i_{\Sigma}^s\bigl(\cl_{\Sigma(s)}(\mathbb R) \cup \cl_{\Sigma(s)}(\mathbb R+i\pi)\bigr)\right) \cup \partial \mathbb D \setminus S$.

The proof of the theorem is complete.

\subsection{Proofs of Theorems \ref{cechthm}, \ref{concompthm} and Corollary \ref{finalcor}}

\begin{proof}[Proof of Theorem \ref{cechthm}]
(1) Consider first the case of $S$ being a finite subset of $\partial \mathbb D$. For $s\in S$ we define 
\begin{equation*}
U_1:=b_\Sigma^S(\mathbb D) \setminus (a_\Sigma^S)^{-1}(s).
\end{equation*}
Let $U_2$ be the union of $(a_\Sigma^S)^{-1}(s)$ and a circular neighbourhood of $s$ whose closure is a proper subset of $\bar{\mathbb D}$.
Both $U_1$, $U_2$ are open in $b_\Sigma^S(\mathbb D)$ and $U_1 \cap U_2=U_2\setminus (a_\Sigma^S)^{-1}(s)$ is the circular neighbourhood of $s$. Since $U_1 \cap U_2$ is contractible, one has $H^k(U_1 \cap U_2,\mathbb Z)=0$, $k \geq 1$. Let us show that for any $k\in\mathbb Z$,
\begin{equation}
\label{reqiso}
H^k(U_2,\mathbb Z) \cong H^k(b_{\Sigma(s)}(T),\mathbb Z).
\end{equation}
To this end consider a sequence $V_1\supset V_2\supset \dots $ of circular neigbourhoods of $s$ such that $\cap_{k=1}^\infty \bar{V}_k=\{s\}$
and $V_1=U_1\cap U_2$. We set $$\hat U_k:=V_k\cup (a_\Sigma^S)^{-1}(s).$$
Let $\iota_{l}^m:\hat U_m\hookrightarrow\hat U_l$, $m\ge l$, be the corresponding embedding. Then $(a_\Sigma^S)^{-1}(s)$ is the inverse limit of the family $\{\hat U_j\, ;\, \iota\}_{j\in\mathbb N}$. It is well known (see, e.g., \cite{Bre}, Chapter II, Corollary
14.6) that the direct limit of \v{C}ech cohomology groups $ H^k(\hat U_l,\mathbb Z)$ with respect to this family gives   $H^{k}((a_\Sigma^S)^{-1}(s),\mathbb Z)$. Note also that
each $\hat U_l$ is a deformation retract of $\hat U_1:=U_2$. Thus the maps $\iota_1^l$ induce isomorphisms $H^k(U_2,\mathbb Z)\cong H^k(\hat U_l,\mathbb Z)$, $l\in\mathbb N$. Since $(a_\Sigma^S)^{-1}(s)\cong b_{\Sigma(s)}(T)$, these facts imply \eqref{reqiso}.

Further, consider the Mayer-Vietoris sequence corresponding to cover $\{U_1,U_2\}$ of $b_{\Sigma}^S(\mathbb D)$:
\begin{equation*}
\dots \to H^{k-1}(b_\Sigma^S(\mathbb D),\mathbb Z) \to H^k(U_1 \cap U_2,\mathbb Z) \to H^k(U_1,\mathbb Z) \oplus H^k(U_2,\mathbb Z) \to H^k(b_\Sigma^S(\mathbb D),\mathbb Z) \to \dots .
\end{equation*}
By the above results $H^k(U_1 \cap U_2,\mathbb Z)=0$ and $H^k(U_2,\mathbb Z)\cong H^k(b_{\Sigma(s)}T,\mathbb Z)$, $k \geq 1$.
Therefore, 
$$H^k(b_\Sigma^S(\mathbb D),\mathbb Z)=H^k(U_1,\mathbb Z) \oplus H^k\bigl(b_{\Sigma(s)}(T),\mathbb Z\bigr),\quad k\ge 1.$$ Proceeding further inductively (i.e., applying similar arguments to $U_1$ etc.) and using the fact that $H^k(b_\Sigma^S(\mathbb D)\setminus S,\mathbb Z)=0$, $k \geq 1$, we obtain that $$H^k(b_\Sigma^S(\mathbb D),\mathbb Z)=\bigoplus_{s \in S} H^k\bigl(b_{\Sigma(s)}(T),\mathbb Z\bigr).$$

Now, if $S \subset \partial \mathbb D$ is an arbitrary closed subset, then since $b_\Sigma^S(\mathbb D)$ is the inverse limit of
$b_{\Sigma|_{F}}^F(\mathbb D)$ for all possible finite subsets $F\subset S$, by the cited result in \cite{Bre} $H^k(b_\Sigma^S(\mathbb D),\mathbb Z)$ is the direct limit of $H^k(b_{\Sigma|_{F}}^F(\mathbb D),\mathbb Z)$. Based on the case considered above we obtain that this limit is isomorphic to $\bigoplus_{s \in S} H^k\bigl(b_{\Sigma(s)}(T),\mathbb Z\bigr)$.

This proves the first part of the theorem.

(2) As is shown in \cite{Br}, if $\Gamma \subset \mathbb R_+$ or $\Gamma \subset \mathbb R_-$, then $b_\Gamma(T)$ is contractible. Therefore under hypotheses of the theorem $H^k\bigl(b_{\Sigma(s)}(T),\mathbb Z\bigr)=0$ for all $s \in S$. The required result now follows from (1),
i.e., $H^k(b_\Sigma^S(\mathbb D),\mathbb Z)=0$ for all $k\ge 1$.

Further, according to \cite{BS} the connectedness of $b^S_\Sigma(\mathbb D)$ and the topological triviality of any complex vector bundle of a finite rank over $b^S_\Sigma(\mathbb D)$ are sufficient for projective freeness of $SAP_\Sigma(S) \cap H^\infty(\mathbb D)$. 

Clearly $b^S_\Sigma(\mathbb D)$ is connected. For otherwise, according to the Shilov theorem on idempotents, see \cite{Sh}, $SAP_\Sigma(S)\cap H^\infty(\mathbb D)$ contains a function $f$ not equal identically to $0$ or $1$ on $\mathbb D$ such that $f^2=f$, a contradiction. 

Next, we show that any finite rank complex vector bundle $\xi$ over $b^S_\Sigma(\mathbb D)$ is topologically trivial. 

Since $b_\Sigma^S(\mathbb D)$ is the inverse limit of the system $\{b_{\Sigma|_{F}}^F(\mathbb D)\, ;\, \omega\}_{F\subset S\, ;\, \#F<\infty}$, see Section~3.4, $\xi$ is isomorphic (as a topological bundle) to pullback to $b_\Sigma^S(\mathbb D)$ of a bundle on some $b_{\Sigma|_{F}}^F(\mathbb D)$ with $F\subset S$ finite, see, e.g., \cite{ES} and \cite{Hu}.
Therefore it suffices to prove the statement for $S\subset\partial\mathbb D$ being a finite subset. In this case, for each $s\in S$ by the contractibility of $(a_\Sigma^S)^{-1}(s)\cong b_{\Sigma(s)}(T)$ (see \cite{Br}) we have that the restriction of $\xi$ to $(a_\Sigma^S)^{-1}(s)$ is topologically trivial. Using a finite open cover $\{U_i\}_{1\le i\le m}$ of $(a_\Sigma^S)^{-1}(s)$ such that $\xi|_{U_i}\cong U_i\times\mathbb C^n$, $n=\rank_{\mathbb C}\,\xi$, for each $i$, we extend (by the Urysohn lemma) global continuous sections $t_{j}: (a_\Sigma^S)^{-1}(s)\to\xi$, $1\le j\le n$, determining the trivialization of $\xi$ over $(a_\Sigma^S)^{-1}(s)$ to each $U_i$. Then using a continuous partition of unity subordinate to a finite refinement of $\{U_i\}_{1\le i\le m}$ we glue together these extensions to get global continuous sections $\tilde t_j$, $1\le j\le n$, of $\xi$ on a neighbourhood $U_s$ of $(a_\Sigma^S)^{-1}(s)$ in $b_\Sigma^S(\mathbb D)$ such that $\tilde t_j|_{(a_\Sigma^S)^{-1}(s)}=t_j$ for each $j$. Since sections $t_j$, $1\le j\le n$, are linearly independent at each point of $(a_\Sigma^S)^{-1}(s)$, diminishing, if necessary, $U_s$ we obtain that sections $\tilde t_j$, $1\le j\le n$, are linearly independent at each point of $U_s$. Thus $\xi$ is topologically trivial on $U_s$. Also, by the definition of the topology on $b_\Sigma^S(\mathbb D)$ without loss of generality we may assume that $U_s\setminus (a_\Sigma^S)^{-1}(s)$ is a circular neighbourhood of $s$.

Suppose that $S=\{s_1,\dots, s_k\}$. Let us cover $b_\Sigma^S(\mathbb D)$ by sets $U_j:=U_{s_j}$, $1\le j\le k$, described above and by 
$U_0:=\bar{\mathbb D}\setminus V$, where $V\subset\cup_{j=1}^k U_{s_{j}}$ and $V\cap U_{s_j}$ is a circular neighbourhood of $s_{j}$ distinct from $U_{s_j}\setminus (a_\Sigma^S)^{-1}(s_j)$, $1\le j\le k$. Since $U_0$ is contractible, $\xi|_{U_0}$ is topologically trivial. Using trivializations of $\xi$ on $U_{j}$, $0\le j\le k$, we obtain that $\xi$ is defined by a $1$-cocycle $\{c_{ij}\}$ with values in $GL_n(\mathbb C)$ defined on intersections $U_i\cap U_j$, $0\le i<j\le k$. 
In turn, by the definition of sets $U_j$, there is an acyclic cover $\{\tilde U_j\}_{j=0}^{k}$ of $\bar{\mathbb D}$ such that
$(a_\Sigma^S)^{-1}(\tilde U_j)=U_j$, $0\le j\le k$.
Thus there exists a cocycle $\{\tilde c_{ij}\}$ on $\{\tilde U_j\}_{j=0}^k$ such that $\tilde c_{ij}\circ a_\Sigma^S=c_{ij}$ for all $i,j$.
This cocycle determines a continuous vector bundle $\tilde\xi$ on $\bar{\mathbb D}$ trivial on each $\tilde U_i$, $0\le i\le k$, such that
$(a_\Sigma^S)^*\tilde\xi=\xi$. Since $\bar{\mathbb D}$ is contractible, $\tilde\xi$ is topologically trivial. Hence $\xi$ is topologically trivial as well.

The proof of the theorem is complete.
\end{proof}

\begin{proof}[Proof of Corollary \ref{finalcor}]
Let $G \subset \mathbb R$ be an additive subgroup.
We denote by $APC_{G}(T) \subset APC(T)$ 
the algebra of uniformly continuous almost periodic functions on $T$ having their spectrum in $G$. Here the spectrum of a function in $APC(T)$ is the union of the spectra of its restrictions to each horizontal line in $T$ (see \cite{Bes}). The vector space of functions
$\sum_{j=1}^k c_j(y)e^{i\lambda_j x}$, $x+iy\in T$, $c_j\in C([0,\pi])$, $\lambda_j\in G$, $k\in\mathbb N$, is dense in $APC_G(T)$ and, hence, the maximal ideal space $M(APC_{G}(T))$ of $APC_{G}(T)$ is homeomorphic to $b_G(\mathbb R)\times [0,\pi]$. On the other hand, $APH_G(T)\subset APC_G(T)$ and the extension of $APH_G(T)$ to $M(APC_{G}(T))$ separates the points of $M(APC_{G}(T))$.  Since the image of $T$ in $b_G(T)$ is dense (see the proof of Theorem \ref{coronathm}), the latter implies that $b_{G}(T) \cong M(APC_{G}(T))$.
Hence, taking $G:=\Sigma(s)$, $s\in S$, we obtain $$H^k(b_{\Sigma(s)}(T),\mathbb Z) \cong H^k(b_{\Sigma(s)}(\mathbb R),\mathbb Z).$$ Since $b_{\Sigma(s)}(\mathbb R)$ is a compact connected abelian group, the required statements follow from the remark before the formulation of the corollary, and Theorem \ref{cechthm} (1).
\end{proof}

\begin{proof}[Proof of Theorem \ref{concompthm}]
In what follows we assume that uniform algebras are defined on their maximal ideal spaces via the Gelfand transforms.

We will require the following auxiliary result.

\begin{lemma}\label{L4.12}
Assume that a set-valued map $\Sigma$ as in Section~3.1 is defined on $\{-s,s\}$ and $f\in  SAP_{\Sigma}(\{-s,s\})\cap H^\infty(\mathbb D)$.
Consider the function
$$
H_s f(z):=f\left(\frac{z+s}{2}\right), \quad z \in \mathbb D.
$$
Then $H_s f\in SAP_{\Sigma|_{\{s\}}}(\{s\})\cap H^\infty(\mathbb D)$ and 
\begin{equation}\label{shift}
[(i_{\Sigma|_{\{s\}}}^s\circ\iota_{\Sigma(s)})^{*}H_s f](z)=[(i_{\Sigma}^s\circ\iota_{\Sigma(s)})^{*}f](z-\ln 2),\quad z\in T,
\end{equation}
see Theorem \ref{maxidthm}.
\end{lemma}
This result states that $H_s: SAP_{\Sigma}(\{-s,s\})\cap H^\infty(\mathbb D)\to SAP_{\Sigma|_{\{s\}}}(\{s\})\cap H^\infty(\mathbb D)$ is a bounded linear operator which induces under the identification of the fibre $(a^{\{-s,s\}}_\Sigma)^{-1}(s)$ with $b_{\Sigma(s)}(T)$ by $i_{\Sigma}^s$ the map $APH_{\Sigma(s)}(T)\to APH_{\Sigma(s)}(T)$ defined by $h(z)\mapsto h(z-\ln 2)$, $z\in T$, $h\in APH_{\Sigma(s)}(T)$.
\begin{proof}
Clearly $H_s f$ is holomorphic on $\mathbb D$ and continuous on $\partial\mathbb D\setminus\{s\}$.
Let us consider the function $g(z):=[(H_s f)\circ (\Log\circ\varphi_s)^{-1}](z)- [f\circ (\Log\circ\varphi_s)^{-1}](z-\ln 2)$, $z\in T$. Next, we have 
$$
\frac{(\Log\circ\varphi_s)^{-1}(z)+s}{2}-(\Log\circ\varphi_s)^{-1}(z-\ln 2)=\frac{se^{2z}}{(2i+e^z)(4i+e^z)}\to 0\quad \text{as}\quad
\Real(z)\to -\infty .
$$
Since by the definition of $SAP_{\Sigma}(\{-s,s\})\cap H^\infty(\mathbb D)$ the function $f\circ (\Log\circ\varphi_s)^{-1}$ is uniformly continuous on $T$, from the last expression we obtain that $g(z)\to 0$ as $\Real(z)\to-\infty$. But $\Real(z)\to-\infty$ if and only if
$(\Log\circ\varphi_s)^{-1}(z)\to s$. Therefore the function $g\circ\Log\circ\varphi$ is continuous in a circular neighbourhood of $s$ and equals $0$ at $s$. Since the pullback of the function $[f\circ (\Log\circ\varphi_s)^{-1}](z-\ln 2)$, $z\in T$, by $(\Log\circ\varphi_s)^{-1}$ belongs to $SAP_{\Sigma}(\{-s,s\})\cap H^\infty(\mathbb D)$ (it is obtained as the composition of $f$ with a M\"{o}bius transformation preserving points $-s$ and $s$), the function $H_sf\in SAP_{\Sigma|_{\{s\}}}(\{s\})\cap H^\infty(\mathbb D)$. Now, the identity \eqref{shift} follows from the fact that $(g\circ\Log\circ\varphi)(s)=0$ by the definition of $i_\Sigma^s$.
\end{proof}
\begin{corollary}\label{C4.13}
Let $f\in G_{\Sigma(s)}^n(T)$, see Section~3.5. Consider the function
\[
F:=H_s[Kf(\Log\circ\varphi_s)],\quad \text{where}\quad Kf(z):=f(z-\ln 2),\quad z\in T.
\]
Then $F\in G_{\Sigma|_{\{s\}}}^n(\{s\})$ and
\[
(i_{\Sigma|_{\{s\}}}^s\circ\iota_{\Sigma(s)})^{*}F=f.
\]
\end{corollary}
\begin{proof}
The fact that $F\in G_{\Sigma|_{\{s\}}}^n(\{s\})$ follows from the proof of Lemma \ref{L4.12} because the pullback by $\Log\circ\varphi_s$ maps
$APH_{\Sigma(s)}(T)$ isometrically into $SAP(\{-s,s\})\cap H^\infty(\mathbb D)$ so that $\speca_s$ of each of the pulled back function is a subset of $\Sigma(s)$. The second statement of the corollary follows directly from \eqref{shift} because $(i_{\Sigma|_{\{s\}}}^s\circ\iota_{\Sigma(s)})^*(h\circ\Log\circ\varphi_s)=h$ for any $h\in APH_{\Sigma(s)}(T)$ by Theorem \ref{maxidthm}\, (1).
\end{proof}

We are ready to prove the theorem. First we will consider the case $S=\{s_1,\dots, s_m\}$ a finite subset of $\partial\mathbb D$. 

By the definition of connected components of $GL_n(A)$, where $A$ is a Banach algebra, the map $f\mapsto ((i_{\Sigma|_{\{s_1\}}}^{s_1}\circ\iota_{\Sigma(s_1)})^*f,\dots, (i_{\Sigma|_{\{s_m\}}}^{s_m}\circ\iota_{\Sigma(s_m)})^*f)$, $f\in G_\Sigma^n(S)$, induces a homomorphism 
\[
\Psi_S: [G_\Sigma^n(S)]\to\bigoplus_{s_i\in S}[G_{\Sigma(s_i)}^n(T)].
\]
We will show that $\Psi_S$ is an isomorphism.

Suppose that $(g_1,\dots,g_m)\in \bigoplus_{s_i\in S}G_{\Sigma(s_i)}^n(T)$ represents an element $[g]\in \bigoplus_{s_i\in S}[G_{\Sigma(s_i)}^n(T)]$.
Then according to Corollary \ref{C4.13} for an element
\[
\tilde g :=H_{s_1}[K(\Log\circ\varphi_{s_1})^*g_1]\cdots H_{s_m}[K(\Log\circ\varphi_{s_m})^*g_m]\in G_\Sigma^n(S) 
\]
and each $l\in\{1,\dots, m\}$ we have
\[
(i_{\Sigma|_{\{s_l\}}}^{s_l}\circ\iota_{\Sigma(s_l)})^*\tilde g=c_{1l}\cdots c_{l-1 l}\cdot g_l\cdot c_{l+1 l}\cdots c_{ml},
\]
where every $c_{jl}$ is an invertible matrix. Since the matrix-function on the right-hand side is homotopic to $g_l$, for the element $[\tilde g]\in [G_\Sigma^n(S)]$ representing $\tilde g$, we obtain $\Psi_S([\tilde g])=[g]$. Hence $\Psi_S$ is a surjection.

To prove that $\Psi_S$ is an injection, we require a modification of the construction of Corollary \ref{C4.13}. So suppose that 
$F_{s_l}=H_{s_l}[K(\Log\circ\varphi_{s_l})^*f]$, where $f\in G_{\Sigma(s_l)}^n(T)$. By the definition,
$F_{s_l}(s_j)$, $j\ne l$, are well-defined invertible matrices. Let $M$ be a matrix-function with entries from $A(\mathbb D)$ such that
$M(s_j)=\Log(F_{s_l}(s_j))$, $j\ne l$, and $M(s_l)=0$. (Here the logarithm of an invertible matrix $c$ is a matrix $\tilde c$ such that $\exp(\tilde c)=c$.) Then we have
\begin{itemize}
\item[(1)]
$\tilde F_{s_l}:=F_{s_l}\cdot\exp(-M)\in G_{\Sigma|_{\{s_l\}}}^n(\{s_l\})$ and satisfies
\[
(i_{\Sigma|_{\{s_l\}}}^{s_l}\circ\iota_{\Sigma(s_l)})^{*}\tilde F_{s_l}=f\quad\text{and}\quad \tilde F_{s_l}(s_j)=I_{n},\quad j\ne l
\]
(here $I_n$ is the unit $n\times n$ matrix);
\item[(2)]
$\tilde F_{s_l}$ is homotopic to $F_{s_l}$.
\end{itemize}
Statement (2) follows from the fact that $\exp(-M)$ clearly belongs to the connected component containing $I_n$.

Now, suppose that $f\in G_\Sigma^n(S)$ is such that every  matrix-function $g_l:=(i_{\Sigma|_{\{s_l\}}}^{s_l}\circ\iota_{\Sigma(s_l)})^*f$, $l\in \{1,\dots, m\}$, belongs to the connected component of $G_{\Sigma(s_l)}^n(T)$ containing the unit matrix $I_n$, (i.e., $[f]\in Ker(\Psi_S)$).
We set
\[
G_{s_l}:=H_{s_l}[K(\Log\circ\varphi_{s_l})^*g_l],\quad G:=\prod_{1\le l\le m}G_{s_l},\quad \tilde G:=\prod_{1\le l\le m}\tilde G_{s_l},
\]
where each $\tilde G_{s_l}$ is constructed from $G_{s_l}$ as $\tilde F_{s_l}$ from $F_{s_l}$. 

According to property (1),
\[
(i_{\Sigma|_{\{s_l\}}}^{s_l}\circ\iota_{\Sigma(s_l)})^*\tilde G=g_l,\quad \text{for}\quad l\in\{1,\dots, m\}.
\]
Moreover, property (2) implies that $\tilde G$ is homotopic to $G$. Observe also that each $G_{s_l}$ is homotopic to $I_n$ (because $g_l$ satisfies this property and so the required homotopy is defined as the image of the homotopy between $g_l$ and $I_n$ under the continuous map
$H_{s_l}\circ K\circ (\Log\circ\varphi_{s_l})^*$) and therefore $G$ and $\widetilde G$ are homotopic to $I_n$. Finally, according to our construction $f\cdot\tilde G^{-1}$ is an invertible matrix with entries from $A(\mathbb D)$. Since $\bar{\mathbb D}$ is contractible, each such a matrix is homotopic to $I_n$. These facts imply that $f$ is homotopic to $I_n$, that is $[f]=1\in [G_\Sigma^n(S)]$, where $[f]$ stands for
the connected component containing $f\in G_\Sigma^n(S)$.

So $\Psi_S$ is an injection which completes the proof of the theorem in the case of a finite $S$.

To prove the result in the general case we require the following lemma.

\begin{lemma}
\label{concomplem}
For every $f \in G_\Sigma^n(S)$ there exists $\bar{f} \in G_{\Sigma|_F}^n(F)$, where $F \subset S$ is finite, such that $\bar{f} \in [f]$.
\end{lemma}
\begin{proof}
Let $M_\Sigma^n(S)$ be the Banach algebra of $n\times n$ matrix-functions with entries in $SAP_{\Sigma}(S) \cap H^\infty(\mathbb D)$ equipped with the norm $\|h\|:=\sup_{z\in\mathbb D}\|h(z)\|_2$, $h\in M_\Sigma^n(S)$, where $\|\cdot\|_2$ is the $\ell_2$ operator norm on the complex vector space $M_n(\mathbb C)$ of $n\times n$ matrices.
According to Corollary \ref{thm4}, $f$ can be approximated in $M_\Sigma^n(S)$ by functions from $M_{\Sigma|_{F}}^n(F)$ for
some finite subsets $F\subset S$. Since the connected component $[f]$ is open (because $G_\Sigma^n(S)\subset M_\Sigma^n(S)$ is open), the latter implies the required statement: there exists $\bar{f} \in G_{\Sigma|_F}^n(F)$, where $F \subset S$ is finite, such that $\bar{f} \in [f]$.
\end{proof}

This lemma implies that $[G_{\Sigma}^n(S)]$ is the direct limit of the family $\{[G_{\Sigma|_{F}}^n(F)];\, F\subset S,\, \#F<\infty\}$. Therefore we can define a homomorphism $\Psi_S: [G_{\Sigma}^n(S)]\to\bigoplus_{s\in S}[G_{\Sigma(s)}^n(T)]$ as the direct limit of homomorphisms $\Psi_F$ described above. Then $\Psi_S$ is an isomorphism because each $\Psi_F$ is an isomorphism on each image.

This proves the first statement of the theorem.

The second statement follows from the fact that if $\Sigma(s) \subset \mathbb R^+$ or $\mathbb R^-$, then the maximal ideal space $b_{\Sigma(s)}(T)$ of Banach algebra $APH_{\Sigma(s)}(T)$ is contractible \cite{Br}. Then the result of Arens \cite{A} implies that $G_{\Sigma(s)}^n(T)$ is connected and therefore $[G_{\Sigma(s)}^n(T)]$ is trivial. From here and the first statement of the theorem we obtain that
$[G_\Sigma^n(S)]$ is trivial, or equivalently, that $G_\Sigma^n(S)$ is connected.
\end{proof}


\begin{thebibliography}{1}

\bibitem{A}
R.~Arens, 
\newblock To what extent does the space of maximal ideals determine the algebra? 
\newblock in ``Function Algebras`` (Birtel, ed.), Scott-Foresman, Chicago, 1966.



\bibitem{Bes}
A.\,S.\, {Besicovich},
\newblock {\em Almost periodic functions}.
\newblock Dover Publications, 1958.




\bibitem{BR}
J.~{Bourgain} and O.~{Reinov},
\newblock On the approximation properties for the space $H^\infty$.
\newblock {Matemathische Nachrichten} {\bf 122} (1983), 19--27.

\bibitem{Bre}
G.~E.~{Bredon}, 
\newblock {\em Sheaf theory}. \newblock Second edition. Graduate Texts in Mathematics {\bf 170},
Springer-Verlag, New York, 1997.

\bibitem{Br}
A.~{Brudnyi},
\newblock Contractibility of maximal ideal spaces of certain
algebras of almost periodic functions.
\newblock {Integral Equations and Operator Theory} {\bf 52} (2005), 595--598.


\bibitem{BK}
A.~{Brudnyi} and D.~{Kinzebulatov},
\newblock On uniform subalgebras of ${L}^\infty$ on the unit circle generated
  by almost periodic functions.
\newblock {Algebra and Analysis} {\bf 19} (2007), 1--33.

\bibitem{BK1}
A.~{Brudnyi} and D.~{Kinzebulatov},
\newblock On algebras of holomorphic functions with semi-almost periodic
boundary values.
\newblock  {C. R. Math. Rep. Acad. Sci. Canada}, in press.


\bibitem{BS}
A.~{Brudnyi} and A.~{Sasane},
\newblock Sufficient conditions for projective freeness of Banach algebras, \newblock {J. Funct. Anal.} {\bf 257} (2009), 4003--4014. 


\bibitem{Cohn}
P.M.~{Cohn}, From Hermite rings to Sylvester domains. \newblock {Proc. Amer. Math. Soc.} {\bf 128} (2000), 1899--1904.

\bibitem{E}
P.~{Enflo}, A counterexample to the approximation property in Banach spaces. \newblock {Acta Math.} {\bf 130} (1973), 309--317.

\bibitem{ES}
S.~Eilenberg and N.~Steenrod,
\newblock {\em Foundations of algebraic topology.} \newblock Princeton University
Press, Princeton, New Jersey, 1952.

\bibitem{F}
R.~{Fox}, Homotopy groups and torus homotopy groups.
\newblock{Ann. Math.} {\bf 49} (1948), 471--510.

\bibitem{Gar}
J.~{Garnett},
\newblock {\em Bounded analytic functions}.
\newblock Academic Press, 1981.

\bibitem{Gra}
H.~{Grauert},
\newblock Holomorphe Funktionen mit Werten in komplexen Lieschen Gruppen. (German)
\newblock {Mathematische Annalen}, {\bf 133} (1957), 450--472.


\bibitem{Gro}
A.~{Grothendieck},
\newblock Products tensoriels toplogiques et espaces nucl\'{e}aires.
\newblock {\em Memoirs Amer. Math. Society} {\bf 16}, 1955.


\bibitem{HM}
K.~{Hofmann} and O.~{Mostert},
\newblock {\em Cohomology Theories for Compact Abelian Groups}. 
\newblock Springer-Verlag, 1973.

\bibitem{Hu}
D.~Husemoller, 
\newblock {\em Fibre bundles}. 
\newblock Springer-Verlag, New York, 1994.

\bibitem{QS}
T.Y.~Lam, 
\newblock {\em Serre's Conjecture. Lecture Notes in Mathematics 635}, \newblock Springer-Verlag, 1978.

\bibitem{Lin}
V.~Ya.~{Lin}, 
\newblock Holomorphic fiberings and multivalued functions of elements of a Banach algebra. 
\newblock {\em Func.
Anal. Appl.}, {\bf 7}, 1973, English translation.





\bibitem{L}
J.~Lindenstrauss,
\newblock Some open problems in Banach space theory.
\newblock {S\'{e}minaire Choquet} {\bf 18} (1975), 1--9.


\bibitem{RS}
L.~Rodman and I. Spitkovsky,
\newblock Almost periodic factorization and corona problem.
\newblock {Indiana Univ. Math. J.} {\bf 47} (1998), 243--281.

\bibitem{Rud}
W.~Rudin,
\newblock {\em Fourier Analysis on Groups}.
\newblock Interscience Publishers, 1962. 


\bibitem{Sar}
D.~Sarason,
\newblock Toeplitz operators with semi-almost periodic kernels.
\newblock {Duke Math J.} {\bf 44} (1977), 357--364.

\bibitem{Sh}
G.~E.~Shilov,
\newblock On decomposition of a commutative normed ring in a direct sums of ideals.
\newblock {Mat. Sb.} {\bf 32(74)}:2 (1953), 353--364. 



\end{thebibliography}
\end{document}